\documentclass[11pt]{article}

\RequirePackage[OT1]{fontenc} \RequirePackage[colorlinks]{hyperref}
\RequirePackage{hypernat}

\usepackage{latexsym,enumerate}
\usepackage{amsmath,amsopn,amstext,amscd,amsfonts,amssymb,
pdfsync,amsbsy}
\usepackage{mathrsfs}
\usepackage{fullpage}
\usepackage{graphicx}
\usepackage[usenames]{color}

\usepackage{bbm}

\def\bB{{\mathbb B}}

\def\b1{{\mathbf{1}}}

\def\bZ{{\mathbb Z}}
\def\bE{{\mathbb E}}
\def\bN{{\mathbb N}}

\def\bR{{\mathbb R}}

\def\bS{{\mathbb S}}

\def\bH{{\mathbb H}}

\def\cA{\mathcal{A}}
\def\cB{\mathcal{B}}
\def\cC{\mathcal{C}}

\def\cF{\mathcal{F}}

\def\cH{\mathcal{H}}

\def\cL{\mathcal{L}}

\def\cS{\mathcal{S}}

\def\cX{\mathcal{X}}

\newcommand{\qed}{\quad\hbox{\vrule width 4pt height 6pt depth 1.5pt}}
\newcommand\eref[1]{(\ref{#1})}
\newtheorem{remark}{Remark}
\newtheorem{proposition}{Proposition}
\newtheorem{definition}{Definition}
\newtheorem{theorem}{Theorem}
\newtheorem{corollary}{Corollary}
\newtheorem{lemma}{Lemma}

\def\eps{{\varepsilon}}
\def\ONE{{\mathbbm 1}}

\def\R{\mathbb{R}}

\newcommand{\Lip}{\operatorname{Lip}}
\newcommand{\diam}{\operatorname{Diam}}
\newcommand{\Ker}{\operatorname{Ker}}
\newcommand{\card}{\operatorname{card}}
\newcommand{\Tr}{\operatorname{Tr}}
\newcommand{\Cte}{\operatorname{const}}
\newcommand{\supp}{\operatorname{supp}}
\newcommand{\Id}{\operatorname{Id}}

\def\AA{A}
\def\DD{{\tilde D}}
\def\XX{{\cX}}
\def\bb{\alpha}

\def\bL{L}
\def\pphi{u}
\def\aaa{\kappa}

\begin{document}

\title{Regularity of Gaussian Processes on Dirichlet spaces}

\author{G\'erard Kerkyacharian,
Shigeyoshi Ogawa,
Pencho Petrushev,
Dominique Picard
}
\date{}
\maketitle

%
%


\begin{abstract}
We are interested in the regularity of centered Gaussian processes $(Z_x( \omega ))_{x\in M}$
indexed by compact metric spaces $(M, \rho)$.
It is shown that the almost everywhere Besov space regularity of such a process is (almost) equivalent to
the Besov regularity of the covariance $K(x,y) = \bE(Z_x Z_y)$ under the assumption that
(i) there is an underlying Dirichlet structure on $M$ which determines the Besov space regularity, and
(ii) the operator $K$ with kernel $K(x, y)$ and the underlying operator $A$ of the Dirichlet structure commute.
As an application of this result we establish the Besov regularity of Gaussian processes indexed by
compact homogeneous spaces and, in particular, by the sphere.
\end{abstract}

\noindent
{\textit{Heat kernel, Gaussian processes, Besov spaces.}}\\
{\textit{MSC 58J35 \and MSC 46E35\and MSC 42C15\and MSC 43A85}}
\setcounter{equation}{0}

\smallskip

\noindent
\smallskip
\section{Introduction}
Gaussian processes have been at the heart of probability theory for very long time. There is a huge literature
about it (see among many others  \cite{Lifshits}, \cite{Led}, \cite{LedTala} \cite{Adler}, \cite{Adler1990} \cite{Rosen}). They also have been playing a key
role in applications for many years and seem to experience an active revival in the recent domains of machine learning (see among others \cite{Rasmussen}, \cite{Seeger}) as well as in Bayesian nonparametric statistics (see for instance \cite{VV1}, \cite{VV2}).

In many areas it is important to develop regularization procedures or sparse representations. Finding adequate regularizations as well as  the quantification of the sparsity play an essential role in the accuracy of the algorithms in statistical theory as well as in Approximation theory. A way to regularize or to improve sparsity which is at the same time genuine and easily explainable is to impose regularity conditions.

The regularity of Gaussian processes has also been for a long time in the essentials of probability theory.
It goes back to Kolmogorov in the 1930s (see among many others \cite{Fernique}, \cite{Talagrand}, \cite{Talagrand3} \cite{Ledoux}, \cite{Li}).

In applications,  an important effort has been put on the construction of Gaussian processes on  manifolds or more general domains, with the two especially challenging  examples of  spaces of matrices and  spaces of graphs to contribute to  the emerging field of signal processing on graphs and extending high-dimensional data analysis to networks and other irregular domains.

Motivated by these aspects we explore in this paper the regularity of Gaussian processes indexed by compact metric domains verifying some conditions in such a way that regularity conditions can  be identified.

In effect, to prove regularity properties, we need a theory of regularity, compatible with the classical examples:
Lipschitz properties and differentiability. At the same time we want to be able to handle more complicated geometries.
For this aspect we borrow the geometrical framework developed in \cite{CKP}, \cite{GKPP}.

Many of the constructions for regularity theorems are based on moments bounds for the increments of the process.
Our approach here is quite different, it utilizes the spectral properties of the covariance operator.
In particular, we use the Littlewood-Paley theory (this point of view was implicitly in \cite{CKR}) to show that the Besov space regularity of the process is (almost) equivalent to
the Besov regularity of the covariance operator. Especially, it is shown that the almost everywhere Besov space regularity of such a process is (almost) equivalent to
the Besov regularity of the covariance $K(x,y) = \bE(Z_x Z_y)$.

 It is also important to notice that unlike many results in the literature, the regularity is expressed using the genuine distance of the domain, not  the distance induced by the covariance.

We illustrate our approach by revisiting the Brownian motion as well as the fractional Brownian motion on the interval. We show the standard Besov regularity of these processes but also prove that they can be associated to a genuine geometry which finally appears in a nontrivial way.

We also illustrate our main result on the more refined case of two points homogeneous spaces and the special case
of the unit sphere $\bS^d$ in $\bR^{d+1}$.

In the two subsequent sections, we recall the background informations about Gaussian processes, the geometrical framework introduced in \cite{CKP}, \cite{GKPP}, and how it allows to develop a smooth functional calculus as well as a description of regularity.
In Section~\ref{sec:main}, we state and prove the main result of the paper.
Section~\ref{sec:BM} details the case of the standard Brownian motion and fractional Brownian motion. In this case, the salient fact is not the regularity result (which is known) but the original geometry corresponding to these processes.
Section~\ref{sec:hom-spaces} deals with positive and negative definite functions on two points homogeneous spaces.
Section~\ref{sec:sphere} establishes the Besov regularity of Gaussian processes indexed by the sphere.

\section{Gaussian processes: Background}\label{sec:background}

In this section we recall some basic facts about Gaussian processes and establish useful notation.

\subsection{General setting for Gaussian processes}\label{subsec:gen-Gauss-process}

Let $(\Omega, \cA, P)$ be a probability space.
A centered Gaussian process on a set $M$ is a family of random variables $Z_x(\omega)$
with $x\in M$ and $\omega \in \Omega$ such that
for all $n \in \bN$, $x_1,\dots,x_n \in M$, and $\alpha_1,\dots,\alpha_n \in \bR$
$$
\sum_{i=1}^n \alpha_i Z_{x_i} \; \hbox{ is a centered Gaussian random variable}.
$$
The covariance function $K(x, y)$ associated to such a process $(Z_x)_{x\in M}$ is defined by
$$
K(x,y):= \bE(Z_x Z_y) \quad\hbox{for} \;\; (x,y) \in M\times M.
$$
It is readily seen that $K(x, y)$ is
real-valued, symmetric, and positive definite, i.e.
$$ K(x,y)= K(y,x) \in \bR, \quad  \hbox{and} $$
$$
\forall n \in \bN, \; \forall \; x_1, \ldots,x_n \in M, \; \forall \alpha_1, \ldots,\alpha_n \in
\bR, \quad  \sum_{i,j\le n } \alpha_i \alpha_j K(x_i,x_j) \ge 0.
$$
Clearly, $K(x, y)$ determines the law of all finite dimensional random variables $(Z_{x_1},\ldots,Z_{x_n})$.

Conversely, if $K(x,y)$ is a real valued, symmetric, and positive definite function on $M\times M$,
 there exists a unique Hilbert space  $\bH$ of functions on $M$ (the associated RKHS),
for which K is a reproducing kernel,
i.e. $f(x)=\langle f,K(x,\bullet) \rangle_{\bH}$, $\forall f \in \bH,\; \forall x\in M$  (see \cite{ARZ}, \cite{SteinChrist}, \cite{CUS}).
Further, if $(\pphi_i)_{i\in I}$ is an orthonormal basis for $\bH$, then the following representation in $\bH$ holds:
$$
K(x,y)= \sum_{i\in I} \pphi_i(x) \pphi_i(y),\quad \forall x,y \in M.
$$
Therefore, if $(B_i(\omega))_{i\in I}$ is a family of independent $N(0,1)$ variables, then
$$ Z_x(\omega):= \sum_{i\in I} \pphi_i(x) B_i(\omega)$$
is a  centered Gaussian process with covariance $ K(x,y)$.
Thus, this is a version of the previous process $Z_x(\omega)$.

\subsection{Gaussian processes with a zest of topology}\label{GP}

We now consider the following more specific setting.
Let $M$ be a compact space and let $\mu$ be a Radon measure on $(M, \cB)$ with support $M$ and
$\cB$ being the Borel sigma algebra on $M$.
Assuming that $(\Omega, \cA, P)$ is a probability space
we let
$$
Z : (M, \cB)\otimes (\Omega, \cA) \mapsto Z_x(\omega) \in \bR, \;
\hbox{be a measurable map}
$$
such that $(Z_x)_{x\in M} $ is a Gaussian process.
In addition, we suppose that $K(x,y)$ is a symmetric, continuous, and positive definite function on $M\times M$.
Then obviously the operator $K$ defined by
$$
Kf(x):= \int_M K(x,y) f(y) d\mu(y), \quad f \in \bL^2(M,\mu),
$$
is a self-adjoint compact positive operator (even trace-class) on $\bL^2(M, \mu)$.
Moreover, $K(\bL^2)\subset C(M)$, the Banach space of continuous functions on $M$.
Let $\nu_1 \ge \nu_2 \ge \dots  > 0$ be the sequence of eigenvalues of $K$
repeated according to their multiplicities
and let $(\pphi_k)_{k\ge 1}$ be the sequence of respective normalized eigenfunctions:
$$
\int_M K(x,y) \pphi_k(y) d\mu(y) = \nu_k \pphi_k(x).
$$
The functions $\pphi_k$ are continuous real-valued functions and the sequence $(\pphi_k)_{k\ge 1}$
is an orthonormal basis for $\bL^2(M,\mu)$.
By Mercer Theorem we have the following representation:
$$
K(x,y)= \sum_k \nu_k \pphi_k(x) \pphi_k(y),
$$
where the convergence is uniform.

Let $\cH \subset \bL^2(\Omega,P)$ be the closed Gaussian space spanned by
finite linear combinations of $(Z_x)_{x\in M}$.
Clearly, interpreting the following integral as Bochner integral with value in the Hilbert space
$\cH$,
we have
$$
B_k(\omega)= \frac 1{\sqrt{\nu_k}} \int_M Z_x(\omega) \pphi_k(x)d\mu(x) \in \cH.
$$
Furthermore,
\begin{align*}
\bE(B_k)&= \bE \Big( \frac 1{\sqrt{\nu_k}} \int_M Z_x(\omega) \pphi_k(x)d\mu(x)\Big)
= \frac 1{\sqrt{\nu_k}} \int_M  \bE (Z_x)   \pphi_k(x) d\mu(x) = 0
\end{align*}
and
\begin{align*}
\bE(B_k B_l)&= \frac 1{\sqrt{\nu_k}\sqrt{\nu_l}}
\bE \Big(\int_M Z_y(\omega) \pphi_k(y) d\mu(y) \int_M Z_x(\omega)  \pphi_l(x)d\mu(x)\Big)\\
&=\frac 1{\sqrt{\nu_k}\sqrt{\nu_l}}  \int_M \int_M \bE (Z_x(\omega)  Z_y(\omega) )  \pphi_k(y) \pphi_l(x)  d\mu(y)  d\mu(x)\\
&=\frac 1{\sqrt{\nu_k}\sqrt{\nu_l}}  \int_M\int_M K(x,y)\pphi_k(x) \pphi_l(y) d\mu(x) d\mu(y)
=\left\{
  \begin{array}{lll}  &1 \;\; \hbox{if}\;\; k=l\\
                      &0 \;\; \hbox{if}\;\; k\ne l.\end{array}
\right.
\end{align*}
As the $B_k$'s belong to the Gaussian space $\cH$,  $B_k$ is a sequence of independent $N(0,1)$ variables.
It is easy to see that
\begin{equation}\label{modif}
\Big\| Z_y- \sum_k \bE (Z_y B_k) B_k\Big\|_{\bL^2({P})}=0 \quad \forall y\in M.
\end{equation}
Indeed, clearly $ \bE(Z_y^2)= K(y,y)$ and
$$
\bE (Z_y B_k)
=\bE \Big(Z_{y}  \frac 1{\sqrt{\lambda_k}} \int_M Z_x(\omega) \pphi_k(x)d\mu(x) \Big)
= \frac 1{\sqrt{\nu_k}}
 \int_M K(y, x) \pphi_k(x) d\mu(x) = \sqrt{\nu_k} \pphi_k(y).
$$
Hence
$$
\bE(Z_y^2)= K(y,y)=\sum_k \nu_k \pphi_k^2(y)=\sum_k ( \bE (Z_y B_k) )^2,
$$
which implies (\ref{modif}).
As a consequence, the process
$$\tilde{Z}_x(\omega) := \sum_k \sqrt{\nu_k} \pphi_k(x) B_k(\omega)$$
is also a modification of $Z_x(\omega)$, i.e.
$P(Z_x= \tilde{Z}_x)=1$, $\forall x \in M$.
\\
We are interested in the regularity of the "trajectory" $x \in M \mapsto Z_x(\omega)$
for almost all $\omega \in \Omega$ and  for a suitable modification of $Z_x(\omega).$
In fact, we  will focus on  the version $ \tilde{Z}_x(\omega) .$

\section{Regularity spaces on metric spaces with Dirichlet structure}\label{sec:setting}

On a compact metric space $(M, \rho)$ one has the scale of $s$-Lipschitz spaces defined by the norm
\begin{equation}\label{def-Lip}
\|f\|_{\Lip_s} := \|f\|_\infty + \sup_{x\ne y}  \frac{|f(x)-f(y)|}{\rho(x,y)^s},\quad 0<s\le 1.
\end{equation}

In Euclidian spaces  a function  can be much more regular than Lipschitz, for instance differentiable at different order,
or belong to some Sobolev space, or even in a more refine way to a Besov space.
For this purpose, we consider
metric measure spaces
with Dirichlet structure.
This setting is  rich enough  to develop
a Littlewood-Paley theory in almost complete analogy with the classical case on $\R^d$, see \cite{CKP,GKPP}.
In particular, it allows to develop Besov spaces $B^s_{pq}$ with all set of indices.
At the same time this framework is sufficiently general to cover a number of interesting cases
as will be shown in what follows.
We next describe the underlying setting in detail.

\subsection{Metric spaces with Dirichlet structure}\label{subsec:background}

We assume that $(M, \mu)$ is a compact, connected measure space, where $\mu$ is a Radon measure with support $M$.
Also, assume that $A$ is a self-adjoint non-negative operator with dense domain $D(A) \subset \bL^2(M, \mu)$.
Let $P_t=e^{-t\AA}$, $t>0$, be the associate self-adjoint semi-group.
Furthermore, we assume that $A$ determines a local and regular Dirichlet structure,
see \cite{CKP} and for details
\cite{FUKU} , \cite{Sturm} ,\cite{Sturm0} ,\cite{Sturm1} ,\cite{Sturm3} ,\cite{BH}, \cite{Grigorian}.
In fact, we assume that $P_t$ is a Markov semi-group ($A$ verifies the Beurling-Deny condition):
$$
0\le f\le 1 \;\;\hbox{and}\;\; f\in \bL^2 \;\;\hbox{imply} \;\; 0\le P_tf \le 1,
\quad \hbox{and also}\;\; P_t\ONE_M=\ONE_M \;\; (\hbox{equivalently}\;\;  A\ONE_M=0).
$$
From this it follows that $P_t$ can be extended as a contraction operator on $\bL^p(M, \mu)$ for $1\le p\le \infty$,
i.e. $\|P_tf\|_p\le \|f\|_p$, and $P_tP_sf=P_{t+s}f$, $t, s>0$.

The next assumption is that there exists a sufficiently rich subspace $\DD \subset D(A)$ (see \cite{BH})
such that $ f \in \DD \Longrightarrow f^2 \in D(A)$.
Then we define a bilinear operator ``square gradiant"
$\Gamma: \DD \times \DD \mapsto \bL^1$ by
$$
\Gamma(f,g):= -\frac 12[A(fg) -fA(g)- gA(f)].
$$
Then $\Gamma(f,f) \geq 0$ and
$
\int_M A(f)g d\mu = \int_M\Gamma(f,g) d\mu \quad (\hbox{Integration by part formula)}.
$

\subsection*{\em Main assumptions:}

\begin{enumerate}
\item
Let
\begin{equation}\label{def-dist}
\rho(x,y):= \sup_{\Gamma(f,f) \le 1} f(x)-f(y)\quad\hbox{for $x, y\in M$.}
\end{equation}
We assume that $\rho$ is a metric on $M$ that generates the original topology on $M$.

\medskip

\item {\em The doubling property:}
Denote
 $B(x,r)= \{y\in M: \rho(x,y)<r\}$.
The assumption is that there exists a constant $d >0$ such that
\begin{equation}\label{DOUB}
\mu(B(x,2r)) \leq 2^d \mu(B(x,r)), \quad \forall x\in M, \;\forall r>0.
\end{equation}
This means that $(M, \rho, \mu)$ is a homogeneous space in the sense of Coifman and Weiss \cite{CW}.
Observe that from (\ref{DOUB}) it follows that
\begin{equation}\label{doubling}
\mu(B(x,\lambda r) ) \le c_0\lambda^d \mu(B(x,r))
\quad\hbox{for $x \in M$, $r>0$, and $\lambda >1, c_0= 2^d.$}
\end{equation}
$d$ is a constant playing the role of a dimension.
\item {\em Poincar\'e inequality:} There exists a constant $c>0$ such that
for all $x\in M$, $r>0$, and $f\in \DD$,
$$
\inf_{\lambda \in \bR} \int_{B(x,r)} (f- \lambda)^2d\mu \le cr^2  \int_{B(x,r)} \Gamma(f,f)d\mu.
$$
\end{enumerate}
As a consequence the associated semi-group $P_t=e^{-t\AA}$, $t>0$, consists of integral operators of continuous
(heat) kernel $p_t(x,y) \geq 0,$ with the following properties:

\smallskip

\noindent
(a) Gaussian localization:
\begin{equation}\label{Gauss-local}
\frac{c_1\exp\{-\frac{\rho^2(x,y)}{c_2t}\}}{\sqrt{\mu(B(x,\sqrt t))\mu(B(y,\sqrt t))}}
\le p_t(x,y)
\le \frac{c_3\exp\{-\frac{\rho^2(x,y)}{c_4t}\}}{\sqrt{\mu(B(x,\sqrt t))\mu(B(y,\sqrt t))}}
\quad\hbox{for} \;\;x,y\in M,\,t>0.
\end{equation}

\noindent
(b) H\"{o}lder continuity: There exists a constant $\aaa>0$ such that
\begin{equation}\label{lip}
\big|  p_t(x,y) - p_t(x,y')  \big|
\le c_1\Big(\frac{\rho(y,y')}{\sqrt t}\Big)^\aaa
\frac{\exp\{-\frac{\rho^2(x,y)}{c_2t} \}}{\sqrt{\mu(B(x,\sqrt t))\mu(B(y, \sqrt t))}}
\end{equation}
for $x, y, y'\in M$ and $t>0$, whenever $\rho(y,y')\le \sqrt{t}$.

\smallskip

\noindent
(c) Markov property:
\begin{equation}\label{hol3}
\int_M p_t(x,y) d\mu(y)= 1
\quad\hbox{for $x\in M$ and $t >0$.}
\end{equation}
Above $c_1, c_2, c_3, c_4>0$ are structural constants.

\begin{remark}\label{rem:setting}
The setting described above is quite general.
This setting covers, in particular, the case of compact Riemannian manifolds.
It naturally contains the cases of the sphere, interval, ball, and simplex with weights.
For more details, see \cite{CKP}.
\end{remark}

\noindent
{\em Notation.}
Throughout we will use the notation $|E|:= \mu(E)$ and
$\ONE_E$ will stand for the characteristic function of $E\subset M$.
Also $\|\cdot\|_p=\|\cdot\|_{\bL^p}:=\|\cdot\|_{L^p(M, \mu)}.$
Positive constants will be denoted by $c, c', c_1, C, C',\dots$ and they may vary at every occurrence.
The notation $a\sim b$ will stand for $c_1\le a/b\le c_2$.
As usual we will denote by $\bN$ the set of all natural numbers and
$\bN_0:=\bN\cup \{0\}$.

\smallskip




Although general the setting described above entails a structure, which in particular
allows to develop a complete Littlewood-Paley theory.
Next, we describe some basic traits of this framework (see \cite{CKP,GKPP}).
For any $t >0$ the operator $P_t:=e^{-t\AA}$ is a Hilbert-Schmidt operator:
\begin{equation}\label{Hilbert-Schmidt}
\| e^{-t\AA}\|_{HS}^2:=\int_M\int_M |P_t(x, y)|^2 d\mu(x) d\mu(y) <\infty.
\end{equation}
The doubling property (\ref{DOUB}) implies that
$M$ being compact is equivalent to $\diam (M) <\infty$ as well as to $\mu(M) <\infty$.
It is also equivalent to
\begin{equation}\label{measure-ball-1}
\int_M \mu(B(y, r))^{-1} d\mu(y) <\infty
\quad\hbox{for all  $r >0$.}
\end{equation}
%
From the compactness of $M$ and the fact that $\AA$ is an essentially self-adjoint non-negative operator it follows that
the spectrum of $\AA$ is discrete and of the form:
$0\le \lambda_1 < \lambda_2 < \cdots$.
Furthermore, the respective eigenspaces $\cH_{\lambda_k}:=\Ker (\AA-\lambda_k \Id)$ are finite dimensional and
$$
\bL^2(M, \mu)= \bigoplus_{k\ge 1} \cH_{\lambda_k}.
$$
Denoting by $P_{\cH_{\lambda_k}}$ the orthogonal projector onto $\cH_{\lambda_k}$ the above means that
$
f=\sum_{k\ge 1}P_{\cH_{\lambda_k}}f
$
in $L^2$ for all $f\in \bL^2(M, \mu)$.
In addition,
\begin{equation}\label{dev-A-Pt}
\AA f = \sum_{k\ge 1} \lambda_k P_{\cH_{\lambda_k}}f, \quad \forall f\in D(\AA),
\quad\hbox{and}\quad
P_t f = \sum_{k\ge 1} e^{-t\lambda_k} P_{\cH_{\lambda_k}}f, \quad \forall f\in \bL^2.
\end{equation}
In general, for a function $g\in L^\infty(\R_+)$ the operator $g(\sqrt \AA)$ is defined by
\begin{equation}\label{def-gsqrtA}
g(\sqrt \AA)f:=\sum_{k\ge 1}g(\sqrt{\lambda_k})P_{\cH_{\lambda_k}}f, \quad \forall f\in\bL^2.
\end{equation}

The spectral spaces $\Sigma_\lambda$, $ \lambda >0$, associated with $\sqrt{\AA}$ are defined by
$$
\Sigma_\lambda := \bigoplus_{\sqrt{\lambda_k}\le \lambda} \cH_{\lambda_k}.
$$
Observe that $\Sigma_\lambda \subset C$ and hence $\Sigma_\lambda \subset \bL^p$ for $1\le p\le\infty$.

From now on we will assume that the eigenvalues $(\lambda_k)_{k\ge 1}$
are enumerated with algebraic multiplicities taken into account,
i.e. if the algebraic multiplicity of $\lambda$ is $m$ then $\lambda$ is repeated $m$ times in the sequence
$0\le \lambda_1 \le \lambda_2\le \cdots$.
We let $(u_k)_{k\ge 1}$ be respective real orthogonal and normalized in $\bL^2$ eigenfunctions of $\AA$,
that is, $Au_k=\lambda_ku_k$.

Let $\Pi_\delta (x,y) := \sum_{\sqrt{\lambda_k} \le \delta^{-1}} u_k(x) u_k(y)$, $\delta >0$,
be the kernel of the orthogonal projector onto $\Sigma_{1/\delta}$.
Then as is shown in \cite[Lemma 3.19]{CKP}
\begin{equation}\label{pdelta}
 \Pi_\delta (x,x)  \sim |B(x,\delta)|^{-1}.
\end{equation}
Further, if $N(\delta,M)$ is the covering number of $M$ (or the cardinality of a maximal $\delta-$net),
then
\begin{equation}\label{covn}
  \dim (\Sigma_{\frac 1{\sqrt t}})\sim  \int_M |B(x, \sqrt t)|^{-1} d\mu(x) \sim N(\sqrt t, M)
\sim\| e^{-t \AA}\|^2_{HS}
\le ct^{-d/2}, \quad t>0.
\end{equation}

A key trait of our setting is that it allows to develop a smooth functional calculus.
In particular, if $g\in C^\infty(\R)$ is even, then the operator $g(t \sqrt \AA)$ defined in (\ref{def-gsqrtA})
is an integral operator with kernel $g(t \sqrt \AA)(x, y)$ having this localization:
For any $\sigma>0$ there exists a constant $c_\sigma>0$ such that
\begin{equation}\label{local-ker}
\big|g(t \sqrt \AA)(x, y)\big|
\le c_\sigma |B(x, t)|^{-1}\big(1+t^{-1}\rho(x, y)\big)^{-\sigma},
\quad\forall x, y\in M.
\end{equation}
Furthermore, $g(t \sqrt \AA)(x, y)$ is H\"older continuous.
An immediate consequence of (\ref{local-ker}) is that the operator $g(t \sqrt \AA)$
is bounded on $\bL^p(M)$:
\begin{equation}\label{bounded-oper}
\|g(t \sqrt \AA)f\|_p \le c\|f\|_p, \quad \forall f\in \bL^p(M), \quad 1\le p\le \infty.
\end{equation}
For more details and proofs, see \cite{CKP,GKPP}.

For discretization (sampling) we will utilize {\em maximal $\delta$-nets}.
Recall that a set $\XX\subset M$ is a maximal $\delta$-net on $M$ ($\delta>0$)
if $\rho(x, y)\ge \delta$ for all $x, y\in\XX$, $x\ne y$, and $\XX$ is maximal with this property.
It is easily seen that a maximal $\delta$-net on $M$ always exists.
Of course, if $\delta>\diam (M)$, then $\XX$ will consists of a single point.
The following useful assertion is part of Theorem~4.2 in \cite{CKP}.

\begin{proposition}\label{prop:d-nets}
There exist a constant $\gamma>0$, depending only on the structural constant of our setting,
such that for any $\lambda >0$ and $\delta:=\gamma/\delta$ there exists a $\delta$-net $\XX$
obeying
\begin{equation}\label{d-net}
2^{-1}\|g\|_\infty \le \max_{\xi\in\XX} |g(\xi)| \le\|g\|_\infty,\quad \forall g\in\Sigma_\lambda.
\end{equation}
\end{proposition}

\subsection{Regularity spaces}\label{subsec:reg-spaces}

%
In the general setting described above, the full scales of Besov and Tribel-Lizorkin spaces are available \cite{CKP,GKPP}.
For the purposes of this study we will utilize mainly Besov spaces.

The Sobolev spaces  $W^k_p=W^k_p(\AA)$, $k\ge 1$, $1\le p\le \infty$, are standardly defined by
\begin{equation}\label{def:Wkp}
W^k_p := \big\{f\in D(\AA^{\frac k2}): \|f\|_{W^k_p} := \|f\|_p + \| \AA^{\frac k2}f\|_p <\infty \big\}.
\end{equation}
Consequently, the Besov space $B^s_{pq}=B^s_{pq}$, $s>0$, $1\le p, q\le\infty$, is standardly defined by interpolation
 as in \cite{Peetre}
\begin{equation}\label{def:Bspq}
B^s_{pq} := \big(\bL^p, W^k_p\big)_{\theta, q}, \quad  \theta:=s/k,
\end{equation}
where $\big(\bL^p, W^k_p\big)_{\theta, q}$ is the real interpolation space between $\bL^p$ and $W^k_p$,
see \cite{CKP}.

\medskip

The following Littlewood-Paley decomposition of functions will play an important role in the sequel.
Suppose $\Phi\in C^\infty(\bR)$ is real-valued, even, and such that
$\supp \Phi\subset [-2, 2]$, $0\le \Phi \le 1$, and $\Phi(\lambda)=1$ for $\lambda\in [0,1]$.
Let $\Psi(\lambda):=\Phi(\lambda)-\Phi(2\lambda)$.
Evidently $\supp \Psi \cap \R_+ \subset [1/2, 2]$.
Set
\begin{equation}\label{def-Psi-j}
\Psi_0:=\Phi \quad\hbox{and}\quad \Psi_j(\lambda):= \Psi(2^{-j}\lambda) \quad \hbox{for} \;j\ge 1.
\end{equation}
It is readily seen that $\Psi_0, \Psi \in C^\infty(\R)$, $\Psi_0, \Psi$ are even,
$\supp \Psi_0\subset [-2, 2]$, $\supp \Psi_j \cap \R_+ \subset [2^{j-1}, 2^{j+1}]$, $j\ge 1$,
and
$
\sum_{j\ge 0}\Psi_j(\lambda) =1
$
for $\lambda\in \R_+$.
Consequently, for any $f\in \bL^p(M, \mu)$, $1\le p\le\infty$, ($\bL^\infty:=C$) one has
\begin{equation}\label{LP-decomp}
f=\sum_{j\ge 0}\Psi_j(\sqrt{\AA})f\quad\hbox{in}\quad \bL^p.
\end{equation}
Note that this decomposition also holds for distributions $f\in\cS'$,
naturally defined in the setting of \S\ref{subsec:background}, see \cite{GKPP}.

The following Littlewood-Paley characterization of Besov spaces uses the functions $\Psi_j$ from above:
Let $s>0$ and $1\le p,q \le \infty$. For a function $f\in \bL^p(M,\mu)$ we have
\begin{equation}\label{char-B-sp}
f \in B^s_{p,q} \Longleftrightarrow
\|\Psi_j(\sqrt{A})f\|_p =\eps_j 2^{-js}, \;\;j\ge 0,\;\;\hbox{with}\;\; \{\eps_j\} \in \ell^q.
\end{equation}
Furthermore, if $f\in B^s_{p,q}$, then $\|f\|_{B^s_{p,q}} \sim \|\{\eps_j\}\|_{\ell^q}$.
We refer the reader to \cite{CKP,GKPP} for proofs and more details on Besov spaces in the setting from \S\ref{subsec:background}.

We next clarify the relationship between $B^s_{\infty,\infty}$ and $\Lip s$. 

\begin{proposition}\label{prop:B-Lip}
$(a)$ For any $0<s\le 1$ we have $\Lip s  \subset B^s_{\infty,\infty}$.

$(b)$ Assuming that $\aaa>0$ is the constant from $(\ref{lip})$, then
$B^s_{\infty,\infty} \subset \Lip s$ for $0<s<\aaa$.
\end{proposition}
This claim follows readily from the results in \cite{CKP,GKPP}.

\begin{remark}
In the most interesting case $\aaa =1$, Proposition~\ref{prop:B-Lip} yields
$\Lip s = B^s_{\infty, \infty}$ for $0<s<1$.
\end{remark}

\section{Main result}\label{sec:main}

We consider a centered Gaussian process $(Z_x)_{x\in M}$ with covariance function
$K(x, y):=\bE(Z_xZ_y)$ as  described in \S\;\ref{GP}, indexed by a metric space $M$ with Dirichlet structure
just as  described in \S\;\ref{subsec:background}.
We will adhere to the assumptions and notation from\; \S \;\ref{subsec:background}.

\subsection{Commutation property}\label{subsec:commute}

We now make the {\em fundamental assumption} that $K$ and $A$ commute in the following sense:

\begin{definition}\label{def:commute}
If $K$ is a bounded operator on a Banach space $\bB$ $(K\in \cL(\bB))$ and
$A$ is an unbounded operator with domain $D(A)\subset \bB$,
we say that $K$ and $A$ commute if $K(D(A))\subset D(A)$ and
$$
KAf=AKf, \quad \forall f\in D(A).
$$
\end{definition}

\begin{remark}\label{prop:commute}
Let $A$ be the infinitesimal generator of a contraction semi-group $P_t$.
Then $K$ and $A$ commute in the sense of Definition~\ref{def:commute} if and only if
$$
KP_t=P_tK, \quad \forall t>0.
$$
We refer the reader to \cite{Davies1}, Theorem~6.1.27.
\end{remark}

\smallskip

We now return to the covariance operator $K$ and the underlying self-adjoint non-negative operator $\AA$ from our setting.
In light of Proposition~\ref{prop:commute} our assumption that $K$ and $A$ commute
implies that they have the same eigenspaces.

Recall that the eigenvalues of $\AA$ are ordered in a sequence
$0=\lambda_1 \le \lambda_2 \le \dots$, where the eigenvalues are repeated according to their multiplicities,
and the respective eigenfunctions $(u_k)_{k\ge 1}$ are real-valued, orthogonal, and normalized in $\bL^2$.
Hence
\begin{equation}\label{orth-decomp}
f=\sum_{k\ge 1}\langle f, u_k\rangle u_k, \quad \forall f\in \bL^2(M, \mu).
\end{equation}
Let $(\nu_k)_{k\ge 1}$ be the eigenvalues of the covariance operator $K$.
Thus we have
\begin{equation}\label{eigenvalues}
\AA u_k= \lambda_k u_k\quad\hbox{and}\quad K u_k= \nu_k u_k, \quad k\ge 1.
\end{equation}

\begin{remark}\label{rem:AK}
As a consequence of the commutation property of $K$ and $\AA$,
the operator $\AA K$ is defined everywhere on $\bL^2(M,\mu)$ and is closed
as $K$ is bounded and $A$ is closed.
Therefore, $\AA K$ is a continuous operator from $\bL^2(M,\mu)$ to $\bL^2(M,\mu)$.
Clearly,
$$ KA f = \sum_{k\ge 1} \langle f,u_k\rangle \lambda_k \nu_k u_k
\quad\hbox{for $f\in \bL^2$ and hence}  \quad
\sup_{k\ge 1}  \lambda_k \nu_k =\|KA\|_{ \cL(\bL^2)}  <\infty.
$$
\end{remark}

\subsection{Main Theorem}

We now come to the main result of this article.

\begin{theorem}\label{MAIN}
Let $(Z_x)_{x\in M}$ be a centered Gaussian process with covariance function $K(x, y):=\bE(Z_xZ_y)$
indexed by a metric space $M$ with Dirichlet structure induced by a self-adjoint operator $\AA$ such that $K$ and $\AA$  commute
in the sense of Definition~\ref{def:commute}.
Then the following assertions hold:

$(a)$ If the covariance kernel $K(x,y)$ has the regularity described by
$$
\sup_{x\in M} \| K(x,\bullet)\|_{B^s_{\infty, \infty}} < \infty
\quad\hbox{for some}\;\; s>0,
$$
then the Gaussian process $Z_x(\omega)$ has the following regularity:
For any $0< \bb <\frac s2$
$$
Z_x(\omega) \in   B^\bb_{\infty,1}
\quad\hbox{for almost all}\; \omega
\quad(B^\bb_{\infty,1}\subset B^\bb_{\infty,\infty}).
$$

$(b)$ Moreover, there exists a unique probability measure  $Q_\alpha$ on the Borelian sets of
$B^\bb_{\infty, 1}$ such that the (canonical)  evaluation process:
$$
\forall x\in M, \; \delta_x: \omega \in B^\bb_{\infty,1} \mapsto   \omega(x)
$$
is a centered Gaussian process of covariance
$$
K(x,y)= \int_{B^\bb_{\infty,\infty}} \delta_x(\omega)\delta_y(\omega) dQ_\alpha(\omega)
\quad\hbox{$[$Ito-Nisio representation$]$}.
$$

$(c)$
Conversely, suppose there exists $\bb>0$ such that
$Z_x(\omega) \in B^\bb_{\infty,\infty}$ for almost all $\omega$.
Then
$$
\sup_{x\in M} \| K(x,\bullet)\|_{B^{2\bb}_{\infty, \infty}} < \infty.
$$
\end{theorem}

\begin{remark}
A key point is that in the above theorem the Besov space smoothness parameter $s>0$ can be arbitrarily large,
while $0<s\le 1$ in the case when the regularity is characterized in terms of Lipschitz spaces.
\end{remark}

For the proof of this theorem we need some preparation.

\subsection{Uniform Besov property of $K(x,y)$ and discretization}

Observe that since the covariance function $K(x,y)$ is a continuous positive definite function on $M\times M$,
then from (\ref{eigenvalues}) it follows that
\begin{equation}\label{rep-K}
K(x,y)= \sum_k \nu_k u_k(x) u_k(y)\quad\hbox{and}\quad \nu_k \ge 0.
\end{equation}
We next represent the Besov norm of $K(x, \bullet)$ in terms of
the eigenvalues and eigenfunctions of $K$ and~$\AA$.

\begin{theorem}\label{Kbes}
Let $s>0$. Then
\begin{equation}\label{K-equiv}
\sup_{x \in M} \| K(x,\bullet) \|_{B^s_{\infty, \infty}}
\sim
\max\Big\{
\sup_{x \in M} \sum_{k: \sqrt{\lambda_k} \le 1}\nu_k u_k^2(x), \;
\sup_{j\ge 1}2^{js}\sup_{x \in M} \sum_{k: 2^{j-1} < \sqrt{\lambda_k} \le 2^{j}}\nu_k u_k^2(x)
\Big\}.
\end{equation}
\end{theorem}

%
%

\noindent
{\bf Proof.}
Note first that from (\ref{char-B-sp}) it follows that (with $\Psi_j$ from (\ref{def-Psi-j}))
$$
\sup_x\|K(x, \bullet)\|_{B^s_{\infty,\infty}}\sim \sup_{j\ge 0} 2^{js}\sup_x \|\Psi_j(\sqrt{A})K(x, \bullet)\|_\infty.
$$
But, using (\ref{rep-K}) we have
$
\big(\Psi_j(\sqrt{A})K(x, \bullet)\big)(y) = \sum_k\Psi_j(\sqrt{\lambda_k})\nu_k u_k(x)u_k(y)
$
and hence, applying the Cauchy-Schwartz inequality it follows that
\begin{align*}
\sup_{x, y}\big|\big(\Psi_j(\sqrt{A})K(x, \bullet)\big)(y)\big|
= \sup_{x}\sum_k\Psi_j(\sqrt{\lambda_k})\nu_k u_k^2(x).
\end{align*}
Consequently,
\begin{align}\label{equiv-1}
\sup_x\|K(x, \bullet)\|_{B^s_{\infty,\infty}}
\sim \sup_j2^{js}\sup_{x}\sum_k\Psi_j(\sqrt{\lambda_k})\nu_k u_k^2(x).
\end{align}
Clearly, from (\ref{def-Psi-j}) we have $0\le \Psi_j\le1$,
$\supp \Psi_0\cap \R_+\subset [0, 2]$, and $\supp \Psi_j\cap \R_+\subset [2^{j-1}, 2^{j+1}]$ for $j\ge 1$.
Therefore,
\begin{align*}
\sup_{x}\sum_k\Psi_0(\sqrt{\lambda_k})\nu_k u_k^2(x)
&\le \sup_{x}\sum_{\sqrt{\lambda_k} < 2}\nu_k u_k^2(x)
\quad\hbox{and}\\
\sup_{x}\sum_k\Psi_j(\sqrt{\lambda_k})\nu_k u_k^2(x)
&\le \sup_{x}\sum_{2^{j-1}< \sqrt{\lambda_k} < 2^{j+1}}\nu_k u_k^2(x),
\quad j\ge 1.
\end{align*}
These estimates and (\ref{equiv-1}) readily imply that the left-hand side quantity in
(\ref{K-equiv}) is dominated by a constant multiple of the right-hand side.

In the other direction, observe that by construction
$\Psi_0(\lambda)=1$ for $\lambda\in [0, 1]$ and
$\Psi_{j-1}(\lambda)+\Psi_j(\lambda)=1$ for $\lambda\in[2^{j-1}, 2^j]$, $j\ge 1$.
Hence
\begin{align*}
\sup_{x}\sum_{\sqrt{\lambda_k} \le 1}\nu_k u_k^2(x)
&\le \sup_{x}\sum_k\Psi_0(\sqrt{\lambda_k})\nu_k u_k^2(x)
\quad\hbox{and}\\
\sup_{x}\sum_{2^{j-1}< \sqrt{\lambda_k} \le 2^{j}}\nu_k u_k^2(x)
&\le \sup_{x}\sum_k\Psi_{j-1}(\sqrt{\lambda_k})\nu_k u_k^2(x)
+\sup_{x}\sum_k\Psi_j(\sqrt{\lambda_k})\nu_k u_k^2(x),
\quad j\ge 1.
\end{align*}
These inequalities and (\ref{equiv-1}) imply that
the right-hand side in (\ref{K-equiv}) is dominated by a constant multiple of the left-hand side.
This completes the proof.
\qed

\medskip

The following corollary is an indication of how the Besov regularity relates with the ``dimension" $d$ of the set $M$,
which appears here through the doubling condition (\ref{DOUB}).

\begin{corollary}\label{dec}
Let $\gamma >d$ and $s=\gamma-d$. Then
$$
\nu_k =O\big(\sqrt{\lambda_k}\big)^{-\gamma}  \Longrightarrow \sup_x \| K(x,\bullet) \|_{B^s_{\infty, \infty}}\le c.
$$
\end{corollary}
{\bf Proof.}
If $\nu_k \le c\big(\sqrt{\lambda_k}\big)^{-\gamma}$, then using (\ref{pdelta}) and
(\ref{doubling}) we get for any $j\ge 1$ and $x\in M$
\begin{align*}
\sum_{k: 2^{j-1} \le \sqrt{\lambda_k}  \leq 2^{j} }\nu_k u_k^2(x)
&\le c2^{-\gamma(j+1)} \sum_{k: 2^{j-1} \le \sqrt{\lambda_k} \le 2^{j} } u_k^2(x)
\le c 2^{-\gamma j} \sum_{k: \sqrt{\lambda_k} \le 2^{j}} u_k^2(x)\\
& = c2^{-\gamma j} \Pi_{2^{j}}(x,x)
\leq c2^{-\gamma j} |B(x,2^{-j})|^{-1}
\le c2^{-j(\gamma-d)}.
\end{align*}
A similar estimate with $j=0$ holds for all $k$ such that $\sqrt{\lambda_k}\le 1$.
Then the corollary follows by Theorem~\ref{Kbes}.
\qed

\begin{remark}\label{rem:est-nu-k}
Observe that
$$
\sup_x\sum_{k: 2^{j-1} \leq \sqrt{\lambda_k}
\le 2^j }\nu_k u_k^2(x)  \le c2^{-js}
\Longrightarrow
\sum_{k: 2^{j-1} \le \sqrt{\lambda_k}  \leq 2^j }\nu_k
=\sum_{k: 2^{j-1} \leq\sqrt{\lambda_k} \leq 2^j } \int_M  \nu_k u_k^2(x) d\mu(x)
\le  c2^{-js} |M|.
$$
\end{remark}

\medskip

We will utilize maximal $\delta$-nets on $M$ along with Proposition~\ref{prop:d-nets} for discretization.
For any $j\ge 0$ we denote by $\XX_j$ the maximal $\delta$-net from Proposition~\ref{prop:d-nets}
with $\delta:=\gamma 2^{-j-1}$ such that
\begin{equation}\label{d-net-j}
2^{-1}\|g\|_\infty \le \max_{\xi\in\XX_j} |g(\xi)| \le \|g\|_\infty, \quad \forall g\in\Sigma_{2^{j+1}}.
\end{equation}
The following claim will be instrumental in the proof of Theorem~\ref{MAIN}.

\begin{proposition}\label{prop:dicrete}
We have
$$
\sup_{x \in M} \sum_{k: \sqrt{\lambda_k} \le 1}\nu_k u_k^2(x)
\sim \max_{\xi \in \XX_0} \sum_{k: \sqrt{\lambda_k} \le 1}\nu_k u_k^2(\xi)
$$
and for any $j\ge 1$
$$
\sup_{x \in M} \sum_{k: 2^{j-1} < \sqrt{\lambda_k} \le 2^{j}}\nu_k u_k^2(x)
\sim
\max_{\xi \in \XX_j} \sum_{k: 2^{j-1} < \sqrt{\lambda_k} \le 2^{j}}\nu_k u_k^2(\xi)
$$
with absolute constants of equivalence.
\end{proposition}

This proposition follows readily from the following

\begin{lemma}\label{DPP}
Let $\XX_j$ be the maximal $\delta$-net from above with $\delta:=\gamma2^{-j}$, $j\ge 0$, and let
$$
H(x,y):=\sum_{\sqrt{\lambda_k} \le 2^{j}} \alpha_k u_k(x) u_k(y),
\quad\hbox{where}\quad \alpha_k \ge 0.
$$
Then
$$
\max_{\xi \in \XX_j} H(\xi,\xi) \le \sup_{x,y \in M} |H(x,y) |\le  4\max_{\xi \in \XX_j} H(\xi,\xi).
$$
\end{lemma}

\noindent
{\bf Proof.} Clearly $H(x,y)$ is a positive definite function and hence
$|H(x,y)| \le \sqrt{H(x,x) H(y,y)}$, implying
\begin{equation}\label{est-H}
\max_{\xi, \eta \in \XX_j}  |H(\xi,\eta) | = \max_{\xi \in \XX_j} H(\xi,\xi).
\end{equation}
Evidently, for any fixed $x\in M$ the function $H(x,y) \in \Sigma_{2^{j}}$ as a function of $y$ and by (\ref{d-net-j})
$$
\sup_{y\in M} |H(x,y)| \le 2\max_{\eta \in \XX_j} |H(x,\eta)|.
$$
Now, using that $H(x,\eta)\in \Sigma_{2^j}$ as a function of $x$, we again apply (\ref{d-net-j})
to obtain
\begin{align*}
\sup_{x,y\in M} |H(x,y)|
&\le 2\sup_{x\in M}\max_{\eta \in \XX_j} |H(x,\xi)|
= 2\max_{\eta \in \XX_j} \sup_{x\in M} |H(x,\eta)|\\
&\le 4\max_{\eta \in \XX_j}\max_{\xi \in \XX_j} |H(\xi,\eta)|
= 4\max_{\xi \in \XX_j} H(\xi,\xi).
\end{align*}
Here for the last equality we used (\ref{est-H}).
This completes the proof. \qed

\subsection{Proof of Parts (a) and (c) of Theorem~\ref{MAIN}}\label{subsec:proof-ac}

(a)
Assume $\sup_{x\in M} \| K(x,\bullet)\|_{B^s_{\infty, \infty}} < \infty$.
Let $(B_k(\omega))_{k\ge 1}$ be a sequence of independent $N(0, 1)$ variables.
Then as alluded in \S\ref{GP}
$$
\tilde Z_x(\omega) := \sum_k \sqrt{\nu_k} u_k(x) B_k(\omega)
$$
is also a version of $Z_x(\omega)$.
Let $\Psi_j$, $j\ge 0$, be the functions from (\ref{def-Psi-j}) and observe that
$f\in B^s_{\infty, 1}$ if and only if
$\|f\|_{B^s_{\infty, 1}}\sim \sum_{j\ge 0}2^{js}\|\Psi_j(\sqrt{A})f\|_\infty <\infty$.
Clearly,
\begin{equation}\label{Psi-Z}
\big(\Psi_j(\sqrt{A})\tilde Z_\bullet(\omega)\big)(x)
= \sum_k \Psi_j(\sqrt{\lambda_k})\sqrt{\nu_k} u_k(x)B_k(\omega).
\end{equation}
For each $x\in M$ this is  a Gaussian variable of variance
$$
\sigma^2_j(x) = \sum_k \Psi^2_j(\sqrt{\lambda_k})\nu_ku_k(x)^2 \le c2^{-js}.
$$
Here we used that $\Psi^2_j(\sqrt{\lambda_k})\le 1$,
the assumption   $\sup_{x\in M} \| K(x,\bullet)\|_{B^s_{\infty,\infty}}<\infty$,
and Theorem~\ref{Kbes}.

For any $\bb>0$ we have
\begin{align*}
\bE \Big(\sum_j 2^{j\bb}\|\Psi_j(\sqrt \AA)\tilde Z_\bullet(\omega)\|_\infty\Big)
&= \sum_j 2^{j\bb}  \bE \big(\|\Psi_j(\sqrt \AA)\tilde Z_\bullet(\omega)\|_\infty\big)\\
&\sim  \sum_j 2^{j\bb}  \bE \big(\sup_{\xi \in \XX_j } |\big(\Psi_j(\sqrt \AA)\tilde Z_\bullet(\omega)\big)(\xi)|\big)\\
&\le c\sum_j 2^{j\bb}2^{-js/2} (1+\log (\card(\XX_j))^{1/2}.
\end{align*}
Above for the equivalence we used (\ref{d-net-j}) and
for the last inequality the following  well known lemma (called the Pisier lemma, see for instance \cite{Ta}, lemma A.3.1):
If $Z_1, \dots, Z_N$ are centered Gaussian variables (with arbitrary variances),
then
$$
\bE\big(\max_{1\le k\le N}|Z_k|\big) \le c(1+\log N)^{1/2} \max_k\big(\bE|Z_k|^2\big)^{1/2}.
$$
By (\ref{covn}), we have $\card(\XX_j) \le c2^{jd}$.
Therefore, if $\bb< \frac s2$, then
$$
\sum_j 2^{j\bb}2^{-js/2}(1+\log (\card(\XX_j))^{1/2}
\le c\sum_j 2^{-j(s/2-\bb)}\big(\log (c2^{jd})\big)^{1/2} <\infty.
$$
Consequently,
$
\bE \Big(\sum_j 2^{j\bb}\|\Psi_j(\sqrt \AA)Z_\bullet(\omega)\|_\infty\Big) <\infty
$
and hence $x\mapsto\tilde Z_x(\omega)\in B^\bb_{\infty, 1}$, $0<\bb<s/2$, $\omega$-a.s.

\medskip

(c)
Suppose now that  $ \omega- a.e. , x  \mapsto Z_x(\omega ) \in B^\bb_{\infty, \infty}$, $\bb>0$.
Then by (\ref{Psi-Z}) and (\ref{char-B-sp}):
$$
\sup_j 2^{j\bb}
\Big\|\sum_k \Psi_j (\sqrt{\lambda_k}) \sqrt{\nu_k} u_k(x) B_k(\omega)\Big\|_\infty <\infty,
\quad\omega-{\rm a.s.}
$$
By (\ref{d-net-j}) this is equivalent to
\begin{equation}\label{BIST}
\sup_j 2^{j\bb}  \max_{\xi \in \XX_j} \Big|\sum_k \Psi_j( \sqrt{\lambda_k}) \sqrt{\nu_k} u_k (\xi) B_k(\omega)\Big| <\infty,
\quad \omega-{\rm a.s.}
\end{equation}
However,
$\{2^{j\bb}  \sum_k \Psi_j( \sqrt{\lambda_k}) \sqrt{\nu_k} u_k (\xi) B_k(\omega) \}_{j\in \bN, \xi \in \XX_j}$
is a countable set of Gaussian centered variables.
The Borell-Ibragimov-Sudakov-Tsirelson theorem (see e.g. \cite{Led}, \S7), in particular, asserts that
if $(G_t)_{t\in T}$ is a centered Gaussian process indexed by a countable parameter set $T$
and $\sup_{t\in T} G_t <\infty$ almost surely, then
$\sup_{t\in T} \bE(G_t^2)<\infty$.
Consequently, (\ref{BIST}) implies
$$
\sup_{j\in \bN, \xi \in \XX_j } \bE\Big( 2^{j\bb}\sum_k \Psi_j(\sqrt{\lambda_k})\sqrt{\nu_k} u_k (\xi) B_k\Big)^2<\infty.
$$
Therefore, there exists a constant $C>0$ such that
$$
\max_{\xi \in \XX_j} \sum_k \Psi_j^2( \sqrt{\lambda_k}) \nu_k u^2_k (\xi)  \le  C2^{-2 j\bb}.
$$
But as before, this yields
$$
\max_{\xi \in \XX_0} \sum_{k: \sqrt{\lambda_k}  \le 1 }\nu_k u_k^2(\xi)
\le \max_{\xi \in \XX_0} \sum_k \Psi_0^2(\sqrt{\lambda_k}) \nu_k u^2_k (\xi)
$$
and, for $j\ge 1$,
\begin{align*}
\max_{\xi \in \XX_j}
\sum_{k, 2^{j-1} \leq \sqrt{\lambda_k}  \le 2^j }\nu_k u_k^2(\xi)
\le 2\max_{\xi \in \XX_j} \sum_k \Psi_{j-1}^2(\sqrt{\lambda_k}) \nu_k u^2_k (\xi)
+ 2\max_{\xi \in \XX_j}\sum_k \Psi_j^2(\sqrt{\lambda_k}) \nu_k u^2_k (\xi)
\le  c2^{-2j\bb}.
\end{align*}
Here we used that $\Psi_{j-1}(\lambda)+\Psi_j(\lambda)=1$ for $\lambda\in[2^{j-1}, 2^j]$,
implying $\Psi_{j-1}^2(\lambda)+\Psi_j^2(\lambda)\ge 1/2$.

Finally, applying Proposition~\ref{prop:dicrete} we conclude from above that
$\sup_{x\in M} \| K(x,\bullet)\|_{B^{2\bb}_{\infty, \infty}} < \infty$. \qed

\subsection{Ito-Nisio theorem and construction of the Wiener measure}

\subsubsection{Wiener measure on $B^s_{\infty,1}$ associated to $K$}


We begin with the following

\begin{lemma}\label{lem:B-norms}
Assume $s>0$ and $1\le p \le \infty$, and let $\Psi_j$, $j\ge 0$, be the functions from $(\ref{def-Psi-j})$.
Then
$$
f \in B^s_{p,1}
\Longleftrightarrow  \sum_{j\ge 0}\|\Psi_j(\sqrt{A})f\|_{B^s_{p, 1}} <\infty
\quad\hbox{and}\quad
\|f\|_{B^s_{p,1}}\sim \sum_{j\ge 0}\|\Psi_j(\sqrt{A})f\|_{B^s_{p, 1}}.
$$
\end{lemma}

\noindent
{\bf Proof.}
From (\ref{LP-decomp}) we have for any $f\in \bL^p$
\begin{equation}\label{Lp-decopm-2}
f=\sum_{j\ge 0}\Psi_j(\sqrt{A})f,\quad \forall f\in \bL^p,
\end{equation}
implying
$
\|f\|_{B^s_{p,1}}\le \sum_{j\ge 0}\|\Psi_j(\sqrt{A})f\|_{B^s_{p, 1}}.
$

\smallskip

For the estimate in the other direction, note that by (\ref{char-B-sp})
$$
\|\Psi_j(\sqrt{A})f\|_{B^s_{p, 1}}
\sim \sum_{\ell\ge 0} 2^{\ell s}\|\Psi_\ell(\sqrt{A})\Psi_j(\sqrt{A})f\|_p.
$$
However, $\supp \Psi_j \cap \R_+ \subset [2^{j-1}, 2^{j+1}]$, $j\ge 1$,
and hence
$\Psi_\ell(\sqrt{A})\Psi_j(\sqrt{A})=0$ if $|\ell-j|>1$.
Therefore,
$$
\|\Psi_j(\sqrt{A})f\|_{B^s_{p, 1}}
\sim \sum_{j-1\le \ell \le j+1} 2^{\ell s}\|\Psi_\ell(\sqrt{A})\Psi_j(\sqrt{A})f\|_p.
$$
On the other hand, by estimate (\ref{bounded-oper}) it follows that
$\|\Psi_j(\sqrt{A})g\|_p \le c\|g\|_p$, $\forall g\in\bL^p$, 
and hence\\
$
\|\Psi_\ell(\sqrt{A})\Psi_j(\sqrt{A})f\|_p \le c\|\Psi_j(\sqrt{A})f\|_p,
$
implying
$$
\|\Psi_j(\sqrt{A})f\|_{B^s_{p, 1}}
\le c2^{js}\|\Psi_j(\sqrt{A})f\|_p
\;\;\Longrightarrow \;\;
\sum_{j\ge 0}\|\Psi_j(\sqrt{A})f\|_{B^s_{p, 1}} \le c \sum_{j\ge 0}2^{js}\|\Psi_j(\sqrt{A})f\|_p
\le c\|f\|_{B^s_{p, 1}}.
$$
The proof is complete. \qed

\smallskip

We now precise Theorem \ref{MAIN}, (a) with the following

\begin{proposition}\label{prop:p5} (Ito-Nisio property.)

Under the hypotheses of Theorem~\ref{MAIN} and with the functions $\Psi_j$, $j\ge 0$, from $(\ref{def-Psi-j})$,
if $\sup_{x\in M}\| K(x,\bullet)\|_{B^s_{\infty,\infty}} <\infty$, then
\begin{equation}\label{equiv-expect}
\bE\Big(\sum_{j\ge 0}\|\Psi_j(\sqrt{A}) Z_\bullet(\omega))\|_{B^\alpha_{\infty,1}}\Big)
\sim \bE\Big(\sum_{j\ge 0} 2^{j\alpha} \|\Psi_j(\sqrt{A}) Z_\bullet(\omega)\|_{\infty}\Big) <\infty,
\end{equation}
the map
$$I: \omega\in \Omega \mapsto \sum_j \psi_j(\sqrt{A})Z_\bullet(\omega)(\cdot) \in B^\alpha_{\infty,1}$$
is measurable, the serie is normally convergent in $ B^\alpha_{\infty,1}$,  and the image probability $Q$ on $B^\alpha_{\infty,1}$ satisfies:
$$
\omega \in B^\alpha_{\infty,1}  \xrightarrow{\delta_x}  \omega(x)
$$
is a centered Gaussian process with covariance $K(x,y)$.
\end{proposition}

\noindent
{\bf Proof.}
The equivalence (\ref{equiv-expect}) follows by the proof of Theorem~\ref{MAIN}, (a) and Lemma~\ref{lem:B-norms}.

As is well known, for any Banach space $B$ with a measure space $(\Omega, \cB)$,
if $G$ is a finite set of indices $b_i \in B$ and $X_i(\omega)$ are real-valued measurable functions,
then
$\omega \mapsto \sum_{i\in G} X_i(\omega) b_i$ is measurable from $\Omega$ to $B$.
Hence,
$$
\omega \in \Omega \mapsto \Psi_j(\sqrt{A}) Z_\bullet(\omega)
= \sum_k \Psi_j(\sqrt{\lambda_k})\sqrt \nu_k u_k(\cdot) B_k(\omega)
\in  B^\alpha_{\infty,1}$$
is measurable.
Consequently, by almost everywhere convergence
$$
I: \omega\in \Omega \mapsto \sum_j \Psi_j(\sqrt{A})Z_\bullet(\omega)(\cdot) \in B^\alpha_{\infty,1}
$$
is also measurable, and $I^{*}({P})=Q$ is a probability measure on the Borelian sigma-algebra,
such that under $Q$ the family of random variables $\delta_x$
$$ \omega \in B^\alpha_{\infty,1}  \xrightarrow{\delta_x}  \omega(x)$$
is a centered Gaussian process with covariance
$K(x,y)=\int_{B^\alpha_{\infty,1}} \omega(x) \omega(y)dQ(\omega)$.
\qed

\smallskip

We next take on the uniqueness of $Q$.

\subsubsection{Gaussian probability on Banach spaces. Proof of Theorem~\ref{MAIN} (b)}

For details in this section we refer to \cite{Bogachev}.

Let $E$ be a Banach space and let $\cB(E)$ be the sigma-algebra of Borel sets on $E$.
Let $E^*$ be its topological dual, and assume $\cF$ is a vector space of real-valued functions defined on $E$, and
$\gamma(\cF,E)$ is the sigma-algebra generated by $\cF$.
If $\cF= \cC_d(E,\bR)$ is the vector space of continuous bounded functions on $E$,
then $\gamma(\cC_d(E,\bR),E)= \cB(E)$ is the Borel sigma-algebra.

As is well known the sigma-algebra $\gamma(E^*,E)$ generated by $E^*$ is $\cB(E)$ if $E$ is separable
(By separation $\cB(E)$) is generated by open balls and by separation and Hahn-Banach theorem
open balls are in $\gamma(E^*,E)$).

\begin{proposition}\label{BKS}
Let $E$ be a separable Banach space.
Let $H $ be a vector subspace of $E^*$,  endowed with the $\sigma(E^*,E)$ topology.
Then
$$
H  \hbox{ is closed }\;  \Longleftrightarrow  \;  H\hbox{  is stable by simple limit}.
$$
\end{proposition}

\noindent
{\bf Proof.}
The implication $\Rightarrow$ is obvious.
We now prove $\Leftarrow$.  As $E$ is a separable Banach space,
this is a consequence of:
$$
\forall R>0, B(0,R)=\{f\in E^*: \|f\|_{E^*} \leq R\} \; \hbox{is metrizable for}\;  \sigma(E^*,E)
$$
and, by Banach-Krein-Smulian theorem, $H$ is  $ \sigma(E^*,E)$-closed
if and only if $\forall R>0$, $B(0,R) \cap H$ is $\sigma(E^*,E)$-closed.
But this is clear, since we only have to verify that for every sequence $(f_n) \subset  B(0,R) \cap H$
such that $\lim_{n \mapsto \infty} f_n(x)= f(x)$, $\forall x\in E$,
we have $f \in  B(0,R) \cap H$, which is what is assumed.
\qed

\begin{corollary}
If $E$ is a separable Banach space and $H$ is a vector subspace of $E^*$,
then $\overline{H}^{\sigma(E^*,E)}$ coincides with
the smallest vector space of functions on $E$, stable by simple limit containing $H$.
Moreover,
$$
\gamma(H,E)= \gamma(\overline{H}^{\sigma(E^*,E)},E).
$$
\end{corollary}

\noindent
{\bf Proof.}
Clearly, as $E^*$ is stable by simple limit (by Banach-Steinhauss theorem),
the smallest   vector space of functions on $E$, stable by simple limit containing $ H$ is contained in $ E^*$.
And if $\gamma (H,E)$  the sigma-algebra generated by $H,$  the vector subspace of $E^*,$
$\gamma(H,E)-$measurable is also stable by simple limit.
\qed

\begin{lemma}
Let $E$ be a separable Banach space.
Let $H$ be a subspace of $E^*$ separating $E$.
Then
$$
\gamma(H,E)= \gamma(E^*,E)=\cB(E).
$$
There is at most one probability measure $P$ on the Borel sets of $E$ such that, under $P$,
$\gamma \in H$ is a centered Gaussian variable with a given covariance $K(\gamma,\gamma')$
on $H$.

Moreover if such a probability exists, then
\begin{enumerate}
\item
$E^*$ is a Gaussian space, and $ \overline{E^*}^{\bL^2(E,P)}$ is the Gaussian  space generated by $H$.

\item
There exists $\alpha >0$ such that
\begin{equation}\label{fernique}
\int_{E} e^{\alpha \| x \|_E^2} dP(x) <\infty.
\end{equation}
\end{enumerate}
\end{lemma}

\noindent
{\bf Proof.}
By the Hahn-Banach theorem
$\overline{H}^{\sigma(E^*,E)}= E^*$
and
$$
\gamma(H,E)= \gamma(E^*,E)=\cB(E).
$$
Now, if $K(\gamma,\gamma')$ is a positive definite function on $H$,
it determines an additive function on the algebra of cylindrical sets
related to $H$:
$$
\big\{x\in E: (\gamma_1(x),\dots,\gamma_n(x) ) \in C\big\}, \; \gamma_i \in H, \; C \; \hbox{ Borelian set of}\; \bR^n.
$$
Now, the sigma-algebra generated by this algebra is the Borelian of $E.$

Assume that such a probability $P$ exists .
Let $\cH = E^* \cap  \overline{H}^{\bL^2(E,P)}$. Clearly $\overline{H}^{\bL^2(E,P)}$
is the Gaussian space generated by $H$, and if $(\gamma_n)_{n\ge 1} \in \cH$ is such that
$\forall x \in E$, $\lim_{n\mapsto \infty}\gamma_n(x)= \gamma(x)$ exists,
then clearly $\gamma \in E^*$ by the Banach-Stheinhauss theorem,
and $\gamma \in \overline{H}^{\bL^2(E,P)}$ since a simple limit of random variables in a
closed Gaussian space belongs to this Gaussian space.
Therefore, $\gamma \in \cH$, which by Proposition~\ref{BKS}
implies that $\cH $ is closed.
But $H \subset \cH$ and $\overline{H}^{\sigma(E^*,E)}= E^*$ leads to $\cH= E^*$.

Finally, (\ref{fernique}) is just the Fernique theorem.
\qed

\begin{corollary}
 Let $M$ be a set and let $E$ be a separable Banach space of real-valued functions on $M$.
 Let
 $$
 \forall t \in M, \;  f \in E  \xrightarrow{\delta_t}   f(t) \in \bR.
 $$
 If $ \delta_t \in E^* $,
 then
 $$
 \cB(E)= \gamma (\cH,E), \; \cH=\{ \sum_{\rm finite}  \alpha_i \delta_{t_i}\}.
 $$
 Let $K(x,y)$ be a definite positive function on $M\times M.$
 There is at most one probability measure $P$ on the Borelian sets of $E$ such that, under $P$,
 $(\delta_x)_{x\in M}$ is a Gaussian process, with covariance
 $K(x,y),$ and $E^*$ is a Gaussian space.

 \end{corollary}

We now come to the main assertion here.

\begin{theorem} [Wiener measure]
In the setting defined above, if $K(x,y)$ a continuous positive definite
function on $M$ such that
$\sup_{x\in M} \| K(x,\bullet) \|_{B^s_{\infty,\infty}} <\infty$
and the associated kernel operator $K$ commutes with $A,$
then there is a unique probability measure $Q$ on the Borelian sets of
$B^\alpha_{\infty,1}$, $\alpha<\frac s2$,
such that the family of random variables:
$$\forall x \in M, \;  \omega \in B^\alpha_{\infty,1}  \xrightarrow{\delta_x}   \omega(x) \in \bR $$
is a centered Gaussian process of covariance $K(x,y)$.
\end{theorem}

This theorem holds due to the previous result and the fact that $B^\alpha_{\infty,1}$ is separable.
It also proves Part (b) of Theorem~\ref{MAIN}.

\subsection{Regularity and reproducing kernel Hilbert spaces}

Let $K(x,y)$, $(x,y)\in M \times M$, be a continuous real-valued positive definite kernel on a compact space $M$.
It is well known that $K$ determines a real Hilbert space $\bH_K$ of functions,
for which the evaluation:
$$
\forall x \in M, \;  \delta_x :  f \in \bH_K^* \mapsto f(x) \; \hbox{is continuous.}
$$
Moreover,
$$
y\mapsto K(x,y)= K_x(y) \in \bH_K, \quad \forall f \in \bH_K, \;
\delta_x(f)= \langle f, K_x\rangle_{\bH_K},
\;\;\hbox{and}  \quad (K_x)_{x\in M} \;  \hbox{is a total set in} \;\bH_K.
$$
The space $\bH_K$ is the completion of ${\rm span}\, \{ K(x,\cdot): x\in M \}$,
more precisely
$$
\bH_K^\circ := \Big\{h(y)=\sum_{i\in F} \alpha_i K(x_i,y): \; \|h\|^2_{\bH}
= \sum_{i, j \in F} \alpha_i \alpha_jK(x_i,x_j)
=  \sum_{j \in F} \alpha_j  h(x_j)\Big\}.
$$
It is well known (see e.g. \cite{CUS}) that
$$
\| h \|^2_{\bH}=0 \;\;\hbox{for}\;\; h \in \bH_K^\circ\; \Longleftrightarrow \;  h(y)=0, \; \forall y \in M.
$$
It is also well known (see \cite{Lus-Pages}) that
$$
K(x,y)= \sum_{i\in I} g_i(x) g_i(y)\; \Longleftrightarrow \;
g_i \in \bH_K, \;\forall i \quad \hbox{and} \quad (g_i)_{i\in I}  \; \hbox{is a tight frame for} \; \bH_K.
$$

In our geometric framework, where the regularity spaces are linked to a
suitable symmetric positive operator, and $K(x,y)$ is compatible with the geometry, we have
$$
K(x,y)= \sum_k \nu_k u_k(x) u_k(y)=  \sum_k [\sqrt{ \nu_k}u_k(x)][\sqrt{ \nu_k}u_k(y)],
$$
where $(u_k)_{k\ge 1}$ is an orthonormal basis for $\bL^2(M,\mu)$ consisting of eigenfunctions of $A$
associated to the eigenvalues $(\lambda_k)_{k\ge 1}$.
Therefore, clearly $(\sqrt{ \nu_k}  u_k )_{k\in \bN, \nu_k \neq 0}$ is a tight frame of $\bH$.

Moreover $ (\delta_x)_{x\in M} \subset \bH_K^*$ is dense in  $\bH_K^*$ in the weak $\sigma (\bH_K^*,\bH_K) $ topology.
In fact, the following theorem holds.

\begin{theorem}

$(a)$ Let $\bN(\nu):= \{k \in \bN, \nu_k \neq 0\}$ and define
$$
\cH= \Big\{f: M \mapsto \bR: f(x)= \sum_{k \in \bN(\nu)} \alpha_k \sqrt{\nu_k} \; u_k(x), \;
(\alpha_k) \in \ell^2 \Big\}
\quad \hbox{with inner product}
$$
$$
\langle f,g\rangle_{\cH}
= \Big\langle\sum_{k\in \bN(\nu)} \alpha_k \sqrt{\nu_k} \, u_k(\cdot),
\sum_{k\in \bN(\nu)} \beta_k \sqrt{\nu_k} u_k(\cdot)\Big\rangle_{\cH}
:= \langle (\alpha_k), (\beta_k) \rangle_{\ell^2(\bN(\nu))}.
$$
Then
$\cH$ is a Hilbert space of continuous functions and
$(\sqrt{\nu_k}u_k)_{k\in \bN(\nu)} $ is an orthonormal basis for $\cH$.
In fact $\bH_K= \cH$.

$(b)$ We have for $s>0$
$$
\bH_K \subseteq B^{\frac s2}_{\infty,\infty}
\;\;\Longleftrightarrow\;\;
\sup_{x\in M} \|K(x, \bullet)\|_{B^s_{\infty, \infty}} <\infty.
$$

$(c)$
Let $\alpha <\frac s2$ and denote by $J:\bH_K \mapsto W=B^{\alpha}_{\infty,1}$ the natural injection, and
$J^*: W^* \mapsto \bH_K^*$.
Then
$$
\overline{Im(J(\bH_K)}= {\rm span} \{ u_k \in  W: \;k \in \bN(\nu)\}.
$$
Under the probability $Q_\alpha$ on $W= B^{\alpha}_{\infty,1}, \; W^*$ is a Gaussian space, and
$$
\int_{W} e^{i \gamma(\omega)} dQ_\alpha(\omega)= e^{- \frac 12\| J^*(\gamma)\|_{\bH_K^*}^2}
\;\; \forall \gamma \in W^*
\quad\;
(\gamma \sim N(0,\| J^*(\gamma)\|_{\bH_K^*}^2).
$$
Moreover, $\overline{W^*}^{\bL^2(W,Q_\alpha)}$ is isometrically isomorphic to $\bH_K^* \sim \bH_K$.

\end{theorem}

\noindent
{\bf Proof.}
(a) We have
$$
\sum_{k\in \bN(\nu)} |\alpha_k| \sqrt{\nu_k} \, |u_k(x)|
\le \Big(\sum_{k \in \bN(\nu)} (|\alpha_k|^2\Big)^{\frac 12} \Big(\sum_{k\in \bN(\nu)}\nu_k |u_k(x)|^2\Big)^{\frac 12}
= \| \alpha \|_{\ell^2(\bN(\nu))}\sqrt{K(x,x)}.
$$
Therefore,
$\cH$ is a space of continuous function,
$\cH$ is a Hilbert space, and
$(\sqrt{\nu_k} u_k)_{k \in \bN(\nu) }$ is an orthonormal basis for $\cH$.
Furthermore,
$$
K_x(y)= K(x,y)= \sum_{k\in \bN(\nu)} \nu_k u_k(x) u_k(y)
= \sum_{k\in \bN(\nu)} [\sqrt{\nu_k}u_k(x)][\sqrt{\nu_k} u_k(y)]\in \cH
$$
and
$$
\langle f, K_x\rangle_{\cH}= f(x), \quad \| K_x \|_{\cH}^2
= \sum_{k\in \bN(\nu)} (\sqrt{\nu_k}  u_k(x))^2= K(x,x).
$$
Also, clearly, $\cH= \bH_K$.

\medskip

(b) Suppose that $\sup_{x\in M} \|K(x, \bullet)\|_{B^s_{\infty, \infty}} <\infty$
and let
$f(x)= \sum_{k \in \bN(\nu)} \alpha_k \sqrt{\nu_k} \, u_k(x)$,
where $(\alpha_k)\in \ell^2$.
Then
$$
\Psi_j(\sqrt{A})f(x)= \sum_{k \in \bN(\nu)} \Psi_j(\sqrt{\lambda_k}) \alpha_k \sqrt{\nu_k} \; u_k(x),
$$
implying, for $j\ge 1$,
\begin{align*}
|\Psi_j(\sqrt{A}) f(x)| &\le  \Big(\sum_{k \in \bN(\nu)} |\alpha_k|^2\Big)^{\frac 12}
\Big(\sum_{k \in \bN(\nu)} |\Psi_j(\sqrt{\lambda_k})|^2 \nu_k|u_k(x)|^2\Big)^{\frac 12}
\\
&\le \|f \|_{\bH_K} \Big(\sum_{k:  2^{j-1} \le \lambda_k \leq 2^{j+1}} \nu_k|u_k(x)|^2\Big)^{\frac 12}
\le  c\|f\|_{\bH_K} 2^{-js/2},
\end{align*}
where for the last inequality we used the assumption and Theorem~\ref{Kbes}.
Similarly
$|\Psi_0(\sqrt{A}) f(x)| \le c\|f\|_{\bH_K}$.
Therefore, in light of (\ref{char-B-sp}),
\begin{equation}\label{hB}
\|f\|_{B^{\frac s2}_{\infty, \infty}} \le c\| f \|_{\bH_K}.
\end{equation}

Suppose now that (\ref{hB}) holds. Then for every sequence $(\alpha_k)\in \ell^2$ with $\|(\alpha_k)\|_{\ell^2} \le 1$
we have
$$
\big|\sum_{k \in \bN(\nu)} \Psi_j(\sqrt{\lambda_k}) \alpha_k \sqrt{\nu_k} \; u_k(x)\Big|
\le c2^{-js/2},
\;\; \forall x\in M,
$$
which by duality implies
$$
\Big(\sum_{k \in \bN(\nu)} |\psi_j(\sqrt{\lambda_k})|^2 \nu_k \; |u_k(x)|^2\Big)^{\frac 12}
\le c2^{-js/2},
\quad j\ge 0.
$$
Just as in the proof of Theorem~\ref{Kbes} we get for $j\ge 1$
\begin{align*}
\sum_{k: 2^{j-1} \le \sqrt{\lambda_k} \le 2^{j}}\nu_k u_k^2(x)
\le  \sum_{k \in \bN(\nu)} |\Psi_{j-1}(\sqrt{\lambda_k})|^2 \nu_k \; |u_k(x)|^2
+ \sum_{k \in \bN(\nu)} |\Psi_j(\sqrt{\lambda_k})|^2 \nu_k \; |u_k(x)|^2
\le c2^{-js}
\end{align*}
and similarly
$\sum_{k: \sqrt{\lambda_k} \le 1}\nu_k u_k^2(x) \le c$.
Consequently,
$\sup_{x\in M} \|K(x, \bullet)\|_{B^s_{\infty, \infty}} <\infty$.

\medskip

(c)
Clearly
$\overline{Im(J(\bH_K)}={\rm span}\, \big\{u_k\in W: k \in \bN(\nu)\big\}$
and from the previous results,
under the probability $Q_\alpha $ on $ W= B^{\alpha}_{\infty,1}, \; W^*$ under $Q_\alpha $ is a Gaussian space.
Let now $ F \subset M$ be finite and $ \sum_{i\in F} \alpha_i \delta_{x_i} \in W^* $.
By construction
$$
\int_W \Big(\sum_{i\in F} \alpha_i \delta_{x_i}(\omega)\Big)^2 dQ_\alpha(\omega)
= \sum_{i, j\in F} \alpha_i \alpha_j K(x_i, x_j)
= \Big\| \sum_{i\in F} \alpha_i \delta_{x_i}\Big\|_{\bH_K^*}^2.
$$
We obtain the result by density of the span of $(\delta_x)_{x\in M}$ in $\bH_K^*$.
\qed

\begin{remark}

Let $f \in \bL^2(M,\mu)$.
Clearly
$$
\tilde{f}(\omega): \omega \in  W= B^\alpha_{\infty,1} \mapsto \int_M f(x) \omega(x) d\mu(x)
$$
belongs to $W^*$.
Hence, under $Q_\alpha, \tilde{f}$ is a Gaussian variable and
\begin{align*}
\bE (\tilde{f})^2
&= \int_{W} \Big(\int_M f(x) \omega(x) d\mu(x)\Big)^2 dQ_\alpha (\omega)
= \int_W \int_M f(x) \omega(x) d\mu(x)\int_M f(y) \omega(y) d\mu(y) dQ_\alpha (\omega)
\\
&=\int_M \int_M f(x)f(y) \Big(\int_{W} \omega(x)\omega(y) dQ_\alpha (\omega)\Big)d\mu(x) d\mu(y)
= \langle Kf, f \rangle_{\bL^2(M,\mu)}.
\end{align*}
Consequently,
$$
\int_W e^{i \tilde{f}(\omega)}dQ_\alpha (\omega)
=e^{-\frac 12 \langle Kf, f \rangle_{\bL^2(M, \mu)}}
\quad\hbox{and}\quad
\| J^*(\tilde{f})\|_{\bH_K^*}^2= \langle Kf, f \rangle_{\bL^2(M,\mu)}.
$$

\end{remark}

\section{The examples of  Brownian Motion and fractional Brownian motion}\label{sec:BM}

Here we illustrate our main result (Theorem~\ref{MAIN}) on the example of the standard Brownian motion.

\subsection{Wiener representation of Brownian motion}

Assuming $M=[0,1]$, consider the kernel
$$
K(x,y)= x\wedge y= \frac{1}{2}(x+y-|x-y|).
$$
It is easy to find the eigenfunctions and eigenvalues of the operator $K$ with kernel $K(x, y)$.
Indeed, we wish to find sufficiently smooth solutions $\phi$ of the problem
\begin{align*}
\int_0^1  x\wedge y   \phi(y) dy= \lambda \phi(x), \; \hbox{where} \;\; \phi(0)=0, \;\;\lambda\ne 0.
\end{align*}
Differentiating both sides of the above identity we obtain
$\int_x^1 \phi(y)= \lambda \phi'(x)$, implying $\phi'(1)=0$.
Another differentiation leads to
$\phi''(x)+\lambda^{-1}\phi=0$.
As a result, we obtain the following eigenfunctions and eigenvalues:
$$
\phi_k(x)= \sqrt 2 \sin \big(k+\frac 12\big)\pi x, \quad  \lambda_k =\frac 1{ (\pi (k+\frac 12))^2}, \;\; k=0,1, \dots .
$$
Therefore,
$$
K(x,y)= \sum_{k\ge 0} \frac 2{((k+ \frac 12)\pi)^2}  \sin \big(k+ \frac 12\big )\pi x \sin\big(k+ \frac 12\big)\pi y.
$$
The associated Gaussian process takes the form
$$
Z_x(\omega) = \sum_{k\ge 0} \frac 1{(k+ \frac 12)\pi} \sqrt 2 \sin\big(k+ \frac 12\big)\pi x\cdot B_k(\omega),
\quad\hbox{where}\;\; B_k \sim N(0,1), \; i.i.d.
$$

A natural Dirichlet space (with Neumann-Dirichlet boundary conditions)
is induced by the operator
$$
\AA f:=-f'', \; D(\AA):= \big\{f\in C^2\,]0,1[\,\cap C^{1}[0,1]: \; f(0)=f'(1)=0 \big\}.
$$
Clearly,
$$
\int_0^1 \AA f(x) f(x) dx = \int_0^1f'^2(x) dx, \quad f\in D(\AA),
$$
and
$$
A\big(\sin \big(k+ \frac 12\big)\pi \bullet\big)(x)
= \big(\big(k+ \frac 12\big)\pi\big)^2 \sin \big(k+ \frac 12\big)\pi x.
$$
Also, the distance on $[0, 1]$ is defined by
$$
\rho(x,y)= \sup_{|f'|\leq 1} |f(x)-f(y)| =|x-y|.
$$
In this setting, the Poincar\'{e} inequality and the doubling property are obvious, and clearly
\begin{equation}\label{besov1}
|K(x,y)-K(x,y')| \le |y-y'|,
\quad \hbox{implying} \quad \sup_{x\in M} \| K(x,\bullet) \|_{B^1_{\infty, \infty } } \le 1.
\end{equation}

So far everything looks fine, unfortunately the Dirichlet space induced here
does not verifies all the conditions described \S\ref{subsec:background}, e.g.
the associated semi-group is not Markovian due to the fact that
the function $\ONE :=\ONE_M$ does not belong to $D(A)$.

In the next subsections we will discuss a useful way to circumvent this problem, in particular,
we will identify a Dirichlet space adapted to the framework of Brownian motion.
This will require careful study of positive and negative definite kernels.

 \subsection{Positive and negative definite functions}\label{subsec:PD-ND}

 For this subsection we refer the reader to \cite{BCR},
 \cite{Schoenberg1},
\cite{Bochner}.
 Recall first the definitions of positive and negative definite functions:

\begin{definition}\label{def-PD}
Given a set $M$, a real-valued function $K(x,y)$ defined on $M\times M$
is said to be {\em positive definite} (P.D.), if
$$
K(x,y)=K(y,x), \quad \hbox{and}\quad
\forall \alpha_1,\dots,\alpha_n \in \bR,\;  \forall x_1,\ldots,x_n\in M,
\quad \sum_{i,j=1}^n\alpha_i\alpha_j K(x_i,x_j)\ge 0.
$$
\end{definition}
As shown in \S\ref{subsec:gen-Gauss-process} the following characterization is valid:
$$
K(x,y) \; \hbox{is P.D.}\quad  \Longleftrightarrow \quad K(x,y)= \bE(Z_xZ_y),
$$
where $(Z_x)_{x\in M}$ is a Gaussian process.

For any $u\in M$ we associate to $K(x,y)$ the following P.D. kernel
$$
K_u(x,y) := K(x,y)+ K(u,u)-K(x,u)-K(y,u)= \bE [(Z_x-Z_u)(Z_y-Z_u)].
$$
Clearly,
$$
K_u \equiv K\Longleftrightarrow K(u,u)=0.
$$

\begin{definition}\label{def:ND}
Given a set $M$, a real-valued function $\psi(x,y)$ defined on $M\times M$ is said to be {\em negative definite} (N.D.),
if
$$
\psi(x,y)=\psi(y,x), \forall x, y\in M, \;\; \psi(x,x)\equiv 0, \quad \hbox{and}
$$
$$
\forall \alpha_1,\ldots,\alpha_n \in \bR\;\;\hbox{s.t.}\;\; \sum_i \lambda_i=0, \;\;
\forall x_1,\dots, x_n\in M, \quad
\sum_{i,j=1}^n \alpha_i\alpha_j \psi(x_i,x_j)\le 0.
$$
\end{definition}

The following characterization is valid (see e.g. \cite[Proposition 3.2]{BCR}):
$$
\psi(x,y) \;\hbox{ is N.D.}\quad  \Longleftrightarrow \quad \psi(x,y)= \bE(Z_x-Z_y)^2,
$$
where $(Z_x)_{x\in M}$ is  a Gaussian process.
Consequently, $\psi(x,y) \ge 0$, $\forall x,y \in M$.

From above it readily follows that $\sqrt{\psi(x,y)}$ verifies the triangular inequality:
\begin{equation}\label{tri}
|\sqrt{\psi(x,y)}- \sqrt{\psi(z,y)}| \le \sqrt{\psi(x,z)}, \quad \forall x,y,z \in M.
\end{equation}

The following proposition can easily be verified.

\begin{proposition}\label{prop:PD-ND}
$(a)$
Let $K(x,y)$ be a P.D. kernel on a set $M$, and set
\begin{equation}\label{def:Psi-K}
\psi_K(x,y):= K(x,x)+ K(y,y)- 2K(x,y).
\end{equation}
Then $\psi_K$ is negative definite. The kernel $\psi_K$ will be termed the N.D. kernel associated to $K$.
In fact, if $K(x,y)= \bE(Z_x Z_y)$, then $\psi_K(x,y)= \bE(Z_x-Z_y)^2$.
Furthermore, $\psi_K \equiv \psi_{K_u}$, $\forall u\in M$.

\medskip

$(b)$
Let $\psi$ be a N.D. kernel, and for any $u\in M$ define
$$
N(u,\psi)(x,y):=\frac 12 [\psi(x,u)+\psi(y,u)-\psi(x,y)].
$$
Thus, if $\psi(x,y)= \bE(Z_x-Z_y)^2$, then $N(u,\psi)(x,y) :=\bE\big[(Z_x-Z_u)(Z_y-Z_u)\big]$.
Then $N(u,\psi)$ is P.D.
Moreover,
$$
N(u,\psi_K)=K_u.
$$
\end{proposition}

The next assertion contains our key idea.

\begin{proposition}\label{prop:p8}
Let $\psi(x,y)$ be a real-valued continuous N.D. function on the compact space $M$,
and set
$$
\tilde{K}(x,y):= \frac 1{2|M|} \int_M [\psi(x,u) + \psi(y,u) -\psi(x,y)] d\mu(u).
$$
Then

$(a)$
$\tilde{K}$ is positive definite.

$(b)$
$\ONE$ is an eigenfunction of the operator $\tilde{K}$ with kernel $\tilde K(x, y)$, that is,
$$
\int_M \tilde{K}(x,y) \ONE(y) d\mu(y)
= \int_M \tilde{K}(x,y) d\mu(y)
= C \ONE,
\quad\hbox{with}\quad C= \frac{1}{2|M|}\int_M\int_M \psi(u,y)d\mu(u) d\mu(y) \geq 0.
$$

$(c)$
$$\exists z \in M \;\;\hbox{s.t.}\;\; \tilde{K}(z,z)=0
\;\Longleftrightarrow\;   \tilde{K}(x,y) \equiv 0
\;\Longleftrightarrow \; \psi(x,y) \equiv 0.
$$
\end{proposition}

\noindent
{\bf Proof.}
Parts (a) and (b) are straightforward.
For the proof of (c) we first observe the obvious implications:
$$
\psi(x,y) \equiv 0 \; \Longrightarrow \; \tilde{K}(x,y) \equiv 0 \;
\Longrightarrow \; \tilde{K}(z,z)=0, \; \forall z\in M.
$$

Now, let $\tilde{K}(z,z)=0$ for some $z\in M$.
Then $\frac 1{2|M|} \int_M [\psi(z,u) + \psi(z,u) -\psi(z,z)] d\mu(u)=0$.
By definition $\psi(z,z)=0$ and hence
$\int_M \psi(z,u)d\mu(u)=0$.
However, $\psi(z,u)$ is continuous and $\psi(z,u)\ge 0$.
Therefore, $\psi(z,u)=0$, $\forall u\in M$.
Now, employing (\ref{tri}) we obtain for $ x,y \in M$
$$
\sqrt{\psi(x,y)} = |\sqrt{\psi(x,y)} -\sqrt{\psi(z,y)}| \leq \sqrt{\psi(x,z)} =0,
$$
and hence $\psi(x,y) \equiv 0$.
This completes the proof. \qed

\begin{remark}\label{exp}
The following useful assertions can be found in e.g. \cite{BCR}, \cite{Schoenberg1}, \cite{Bochner}.

Let $\psi (x, y)$, defined on $M\times M$, obey
$\psi(x,y)= \psi(y,x)$, $\forall x,y \in M$, and $\psi(x,x)\equiv 0$.
Then
\begin{align*}
\psi\;\hbox{is N.D.} \; &\Longleftrightarrow  \forall t>0, \;  e^{-t \psi} \;\hbox{is P.D.}\\
\psi\; \hbox{is N.D.} \; &\Longrightarrow  \forall\; 0<\alpha \leq 1, \; \psi^{\alpha} \;\hbox{is N.D.}\\
\psi\; \hbox{is N.D.} \; &\Longrightarrow \; \log(1+\psi) \;\hbox{is N.D.}
\end{align*}
\end{remark}

The following proposition can easily be verified.

\begin{proposition}\label{prop:p9}

Let $M$ be a compact space, equipped with a Radon measure $\mu$.
Assume that $K(x,y)$ is a~continuous P.D. kernel
and let $\psi:=\psi_K$ be the associated to $K(x, y)$ N.D. kernel,
i.e.
$\psi(x,y):= K(x,x)+ K(y,y)- 2K(x,y)$.
Set
$$
K_u(x,y):= \frac 12 [\psi(x,u) + \psi(y,u) -\psi(x,y)],
$$
and
$$
\tilde{K}(x,y):= \frac 1{2|M|} \int_M [\psi(x,u) + \psi(y,u) -\psi(x,y)]  d\mu(u)
= \frac 1{|M|}\int_M K_u(x,y) d\mu(u).
$$
Denote by $K$ and $\tilde K$ the operators with kernels $K(x, y)$ and $\tilde K(x, y)$.

Then
\begin{equation}\label{def.tildeK}
\tilde{K}(x,y)= K(x,y) + |M|^{-1}\Tr(K) - |M|^{-1} K\ONE(x) - |M|^{-1}K\ONE(y).
\end{equation}
Moreover,
$\psi_{\tilde{K}} =\psi$, $\tilde{K}_u=K_u$, and
$$
\tilde{K}\ONE= C \ONE, \quad
C = Tr(K)-\frac{1}{|M|}\int_M\int_M K(x,y)d\mu(x)d\mu(y)
= \frac{1}{2|M|}\int_M\int_M \psi(u,y)d\mu(u) d\mu(y) \ge 0.
$$
In addition,
\begin{equation}\label{k1}
K= \tilde{K} + \Cte.  \; \Longleftrightarrow \; K\ONE = C'\ONE
\end{equation}
and $\Cte. = -|M|^{-1}(\Tr(K)- 2C')$.
\end{proposition}

\noindent
{\bf Proof.}
From the respective definitions, we infer
\begin{align*}
\psi(x,u) + \psi(y,u) -\psi(x,y)
&= [K(x,x)+ K(u,u)- 2K(x,u)]\\
&+ [K(y,y)+ K(u,u)- 2K(y,u)]-[K(x,x)+ K(y,y)- 2K(x,y)]\\
&=2[K(u,u)- K(x,u)-K(y,u) + K(x,y)]
\end{align*}
and hence
\begin{align*}
\tilde{K}(x,y)
&= K(x,y)+\frac 1{|M|} \int_M [K(u,u)- K(x,u)-K(y,u)]d\mu(u)\\
&= K(x,y) + \frac{1}{|M|}\big(\Tr(K) - K\ONE(x) - K\ONE(y)\big).
\end{align*}
The remaining is a consequence of Proposition \ref{prop:p8}. \qed

\begin{remark}\label{knul}
Observe that if $K(x,y)$, $\psi(x,y)$, and $\tilde K(x,y)$ are as in Proposition~\ref{prop:p9},
then
$$
\exists z \in M , \; \tilde{K}(z,z)=0
\; \Longleftrightarrow \;  \tilde{K}(x,y) \equiv 0
\; \Longleftrightarrow \; \psi(x,y) \equiv 0
\; \Longleftrightarrow \; K(x,y)\equiv \Cte.
$$
Indeed,
clearly we have only to show the implication $\psi(x,y) \equiv 0 \Longrightarrow  K(x,y)\equiv \Cte$.
However,
$$
\psi(x,y) \equiv 0 \Longrightarrow K(x,x)+ K(y,y)= 2K(x,y) \leq 2\sqrt{K(x,x)}\sqrt{K(y,y)},
$$
implying
$(\sqrt{K(x,x)}-\sqrt{K(y,y)})^2 \le 0$,
which leads to
$K(x,x) \equiv \Cte.$ and
$$
K(x,y)= \frac{1}{2}(K(x,x)+ K(y,y)) = \Cte.
$$
\end{remark}

\begin{remark}\label{ktilde}
Assume that we are in the geometrical setting described in \S\ref{subsec:background}, associated to an operator $A$.
Just as in \S\ref{subsec:commute}, suppose $K(x, y)$ is a P.D. kernel such that
the associate operator $K$ commutes with $A$.
From (\ref{kernel-A}) we have $A\ONE=0$.

Moreover, it is easy to see that
\begin{equation}\label{kernel-A}
A \ONE_M=0 \quad \hbox{and}\quad \dim \Ker (A) =1.
\end{equation}
Indeed, the Markov property (\ref{hol3}) yields $A\ONE_M=0$.
To show that $\dim \Ker (A) =1$, assume that $Af=0$, $f\in D(A)$. Then $\Gamma(f,f)=0$.
Assume that $f\ne {\rm constant}$. Then $f(x)\ne f(y)$ for some $x, y\in M$, $x\ne y$.
For $\Gamma(f,f)=0$ we have $\Gamma(a f, a f)=0$ for each $a>0$.
Then by (\ref{def-dist}) $\rho(x,y) \ge a|f(x)-f(f)|$, $\forall a>0$, implying
$\rho(x, y)=\infty$, which is a contradiction because $M$ is connected (see \cite{CKP}).
Therefore, $Af=0$ implies $f= \Cte$. and hence $\dim \Ker (A) =1$.

Hence
$$
AK\ONE =  KA\ONE=0.
$$
However, as $\dim \Ker(A)=1$, necessarily $K\ONE= C\ONE$.
Therefore, $K= \tilde{K} + \Cte$.
\end{remark}

\subsection{Back to Brownian motion}

Assume again that
$M=[0,1]$ and $K(x,y):= x\wedge y =\frac{1}{2}(x+y-|x-y|)$.
We will adhere to the notation introduced in \S\ref{subsec:PD-ND}.

The associated to $K(x, y)$  N.D. kernel $\psi(x, y)=\Psi_K(x, y)$ (see (\ref{def:Psi-K})) take the form
\begin{align*}
\psi(x, y)= K(x,x)+K(y, y)-2K(x, y)
\end{align*}
and the induced P.D. kernel $K_u(x, y)$ becomes
\begin{align*}
K_u(x,y)&= \frac 12 [\psi(x,u) + \psi(y,u) -\psi(x,y)]\\
&= K(x,y)- K(x, u)-K(y, u)+ K(u, u)\\
&= \frac 12[|x-u| + |y-u| -|x-y|], \quad u\in [0,1].
\end{align*}
Thus we arrive at the following P.D. kernel
$$
\tilde{K}(x,y):= \int_0^1 K_u(x,y) du =\frac 14[x^2+ (1-x)^2 + y^2+ (1-y)^2 - 2|x-y|].
$$
Denoting by $\tilde K$ the operator with kernel $\tilde K(x, y)$, we have
$$
\tilde{K}\ONE(x) = \frac 12\Big(\int_0^1 \int_0^1 |y-u|dudy +\int_0^1 |x-u|du -\int_0^1|x-y|dy\Big)
= \frac 12 \int_0^1 \int_0^1 |y-u|dudy =  \frac 16.
$$
Further, using that $\int_0^1 \cos k\pi y \, dy = 0$ for $k\in \bN$ we get
\begin{align*}
\tilde{K}(\cos k\pi \bullet)(x)
=  \frac 12 \int_0^1 \int_0^1 |y-u| \cos k\pi y \, dudy - \frac 12\int_0^1|x-y| \cos k \pi y\, dy
\end{align*}
Integrating by parts we obtain
\begin{align*}
\int_0^1|x-y| \cos k \pi y\,dy
&= |x-y|  \frac{\sin k\pi y}{\pi k}\Big|_{y=0}^1
- \int_0^1\big(-\ONE_{[0,x]}(y)+\ONE_{[x,1]}(y)\big) \frac{\sin k\pi y}{\pi k}\,dy\\
&= \frac{1}{\pi k}\Big(\int_0^x \sin k\pi y \, dy - \int_x^1 \sin k\pi y \,dy \Big)\\
&=-\frac{2\cos k\pi x}{(\pi k)^2}  + \frac{1+ (-1)^k}{(\pi k)^2}.
\end{align*}
By the same token
\begin{align*}
\int_0^1 \int_0^1 |y-u| \cos k\pi y\, dudy
&=\int_0^1 \Big(-\frac{2\cos k\pi u}{(\pi k)^2}  + \frac{1+ (-1)^k}{(\pi k)^2}\Big) du
= \frac{1+ (-1)^k}{(\pi k)^2}.
\end{align*}
Putting the above together we infer
\begin{equation*}
\tilde{K}(\cos k\pi \bullet)(x)= \frac{\cos k\pi x}{(\pi k)^2}, \;\; \forall k\in\bN,
\quad\hbox{and}\quad
\tilde{K}\ONE = 1/6.
\end{equation*}
Observe also that the functions
$\{\ONE\}\cup ( \sqrt 2 \cos k\pi x)_{k \in \bN}$ form an orthonormal basis for $\bL^2(0, 1)$.

\smallskip

Let $H^2(0,1)$ be the space of the functions $f \in \bL^2(0, 1)$ that are two times weakly differentiable
and $f',f'' \in \bL^2(0,1)$.
Consider the operator
$$
Af := -f'', \;\; D(A):= \{f \in H^2([0,1]): f'(0)= f'(1)=0\}.
$$
Clearly,
$$
\int_0^1 (Af) g dx=  \int_0^1 f' g' dx = \int_0^1 f Ag dx
$$
and hence $A$ is positive and symmetric.
In fact, $A$ generates a Dirichlet space,
and also
$$
\cos k\pi x \in D(A) \quad\hbox{and}\quad A(\cos k\pi \bullet) (x)= (\pi k)^2 \cos k\pi x, \;\; k\ge 1.
$$
Let
$H^1[0,1]:= \big\{ f \in \bL^2(0,1): \int_0^1 |f'(u)|^2 du <\infty\big\}$.
This defines a Dirichlet form:
$$
A, D(A)= \Big\{f \in H^1(0,1):  \Big|\int_0^1 f'(x) \phi'(x) dx\Big|\le c\|\phi\|_2,\;  \forall \phi \in H^1(0,1)\Big\}.
$$
Thus
$$
\int_0^1 f'(x) \phi'(x) dx = \int_0^1 Af(x)\phi(x)dx
$$
and the distance is defined by
$$
\rho(x,y)= \sup_{\phi \in H^1: |\phi'|\le 1} \phi(x)-\phi(y)=|x-y|.
$$

The Poincar\'{e} inequality is well known to be true in this case.
So we are now in the setting presented above.
$ K(x,\bullet) $ is uniformly $\Lip 1$.
Therefore, $Z_x$ the centered Gaussian process associated to $K$ is almost surely
$\Lip \alpha$, $\alpha < \frac 12$.

The process $Y_x(\omega)= Z_x(\omega)-Z_0(\omega)$ has the same regularity,
and $\frac{1}{2}(|x| + |y| - |x-y|)$ is the associated kernel.
This is the Brownian motion, with the above regularity.

\subsection{Brownian motion and fractional Brownian motion, through analysis on the circle}

The Laplacian on the torus (for instance on $\bR/2\bZ$) is a typical example of an operator generating
a Dirichlet space with all properties that are required for defining a regularity structure.
If we represent $\bR/2\bZ$ by the arc length parametrisation we have a Dirichlet space associated to:
$$
 Af:=- f'', \quad D(A):= \big\{f\in C^2\,]-1,1[\,\cap C^1[-1,1]: f(-1)=f(1),  f'(-1)=f'(1) \big\},
$$
$$
\int_{-1}^1 Af(x) g(x)dx = \int_{-1}^1 f'(x) g'(x) dx,
$$
and the distance is defined by
$$
\rho(x,y)= \inf\big\{f(x)-f(y): |f'|\le 1, f(-1)=f(1), f'(-1)=f'(1)\big\} = |x-y|\wedge (2-|x-y|).
$$
Clearly, the eigenfunctions of $A$ are $(\cos k \pi x )_{k\in \bN_0}$ and $(\sin k\pi x)_{k \in \bN}$.

\subsubsection{Brownian motion on the circle}

Using the Fourier series expansion, we have,
$$
|x| = \frac 12  - \frac 4{\pi^2}  \sum_{n \in \bN} \frac{ \cos(2n+1)\pi x}{(2n+1)^2},
\quad x \in [-1,1].
$$
Hence,
$$
\rho(x,y)=  |x-y|\wedge (2-|x-y|)
= \frac 12 - \frac 4{\pi^2}\sum_{n \in \bN} \frac{\cos(2n+1)\pi(x-y)}{(2n+1)^2},
$$
implying
\begin{align*}
K(x,y)
&:=\frac 12 - |x-y|\wedge (2-|x-y|)
= \frac 4{\pi^2}  \sum_{n \in \bN} \frac{ \cos(2n+1)\pi(x-y)}{(2n+1)^2}\\
& = \frac 4{\pi} \sum_{n \in \bN} \frac{\cos(2n+1)\pi x\cos(2n+1)\pi y}{(2n+1)^2}
+  \frac 4{\pi}  \sum_{n \in \bN} \frac{\sin(2n+1)\pi x\sin(2n+1)\pi y}{(2n+1)^2}.
\end{align*}
From this it follows that the kernel $K(x,y)$ is P.D. and
$$
\psi_K(x,y)= K(x,x)+ K(y,y)- 2K(x,y)= |x-y|\wedge (2-|x-y|)= \rho(x,y).
$$

Therefore, the Gaussian  process $Z_x(\omega)_{x\in [-1,1]}$  associated to $ K(x,y)$
is a Brownian field with respect to $\rho$. Its regularity is $\Lip \alpha$, $\alpha <\frac 12$
with respect to the metric $\rho$.
We can now restrict to $[0,1]$,  $\rho(x,y)= |x-y|$ $\forall x,y \in [0,1]$.
Thus, considering $W_x := Z_x-Z_0$, restricted to $x\in [0,1]$, we get the classical Brownian motion with
$$
W_0:=0,\quad \bE(W_x-W_y)^2=2 |x-y|
$$
and we again obtain its regularity as a byproduct.

\subsubsection{Fractional Brownian motion on the circle and on $[0,1]$}

Let $0< \alpha <1$.
The Fourier series expansion of $|x|^\alpha$ on $[-1, 1]$ takes the form
\begin{align*}
|x|^\alpha= \frac{1}{\alpha +1} + 2 \sum_{k\ge 1} \cos k \pi x  \int_0^1  u^\alpha \cos k \pi u du.
\end{align*}
Integrating by parts we get
\begin{align*}
\int_0^1 u^\alpha\cos k \pi u\, du
= -\frac{\alpha}{k\pi}  \int_0^1 u^{\alpha-1} \sin k \pi u\, du
= -  \frac{\alpha}{(\pi k)^{\alpha +1}}  \int_0^{k\pi} u^{\alpha-1} \sin u \,du.
\end{align*}
Hence,
$$
\frac{1}{\alpha +1} -|x|^\alpha = 2 \alpha   \sum_{k\geq 1} \gamma_k
\frac{ \cos k \pi x }{(\pi k)^{\alpha +1}},
$$
where
\begin{align*}
\gamma_k =\int_0^{k\pi} u^{\alpha-1} \sin u\, du
&= \sum_{j=0}^{k-1} \int_{j\pi}^{(j+1)\pi}  \frac 1{u^{1-\alpha}} \sin u \,du
= \sum_{j=0}^{k-1} (-1)^j \int_{j\pi}^{(j+1)\pi}  \frac 1{u^{1-\alpha}}|\sin u|du
\\
&= \sum_{j=0}^{k-1} (-1)^j \int_{0}^{\pi}  \frac {\sin u}{ (u+j\pi)^{1-\alpha}}du
=: \sum_{j=0}^{k-1} (-1)^j a_j.
\end{align*}
%
Here $a_0 >a_1 >\cdots \ge 0$ and $\lim_{j\to \infty} a_j =0$.
Hence,
$\gamma :=\lim_{j\to \infty} \gamma_j$ exists and $0<\gamma< \pi^{\alpha+1}/(\alpha+1)$.
Therefore, 
\begin{align*}
K_\alpha (x,y)
&=  \frac{1}{\alpha +1} - (|x-y|\wedge (2-|x-y|)^\alpha =  2 \alpha   \sum_{k\geq 1} \gamma_k
\frac{ \cos k \pi (x-y) }{(\pi k)^{\alpha +1}}\\
&= 2 \alpha   \sum_{k\geq 1} \gamma_k
\frac{ \cos k \pi x \cos k \pi y + \sin k \pi x \sin k \pi y   }{(\pi k)^{\alpha +1}}
\end{align*}
is a P.D. kernel compatible with the Dirichlet structure defined by the Laplacian on the circle $\bR/2\bZ$.
Moreover, as $0<\alpha < 1$
$$
|K_\alpha (x,y)-K_\alpha (x,y')| \le  |\rho(x,y)^\alpha-\rho(x,y')^\alpha| \leq \rho(y,y')^\alpha.
$$
Consequently, the associated  Gaussian process $(Z_x(\omega))_{x\in [-1,1]}$ with covariance function $K_\alpha (x,y)$
is $\Lip \beta$ for $\beta <\frac \alpha 2$ on $[-1,1]$.
If we restrict this process to $x\in[0,1]$ as $\rho(x,y)=|x-y|$, $\forall x,y \in [0,1]$,
we get a Gaussian process on $[0,1]$ with covariance $\frac{1}{\alpha +1} - |x-y|^\alpha$, and such that
$$
\bE(Z_x-Z_y)^2= |x-y|^\alpha.
$$
Hence, the process $(Z_x-Z_0)_{x\in [0,1]}$ has covariance $|x|^\alpha +|y|^\alpha-|x-y|^\alpha$
and regularity $\Lip \beta$, $\beta <\frac \alpha 2$.
This is the standard fractional Brownian function.

\begin{remark}\label{rem:frac-Brown}
If $\alpha >1$, then as above
$\gamma_k = \sum_{0\leq j \leq k-1} (-1)^j\int_{0}^{\pi}(u+j\pi)^{\alpha-1} \sin u du$
and hence $\gamma_k >0$ if $k$ is even and $\gamma_k <0$ if $k$ is odd.
From this one can deduce that $\rho(x,y)^\alpha$ is not a definite negative function on the circle.
\end{remark}

\section{Positive and negative definite functions on compact homogeneous spaces}\label{sec:hom-spaces}

Here we present some basic facts about positive and negative definite kernels in
the general setting of compact two point homogeneous spaces.
Then, in the next section, we utilize these results and our main Theorem~\ref{MAIN}
to establish the Besov regularity of Gaussian processes indexed by the sphere.

\subsection{Group acting on a space}

Let $(M, \mu)$ be a compact space equipped with a positive Radon measure $\mu$.
Assume that there exists a~group $G$ acting transitively on $(M,\mu)$, that is,
there exists a map $(g,x) \in G \times M \mapsto  g\cdot x \in M$ such that

\smallskip

1. $h\cdot (g\cdot x)= (hg)\cdot x$, $\forall g,h \in G$,

\medskip

2. $\exists e\in G$ \; s.t. \; $e\cdot x=x$, \; $\forall x\in M$ \; ($e$ is the neutral element in $G$),

\medskip

3. $\forall x,y \in M$, $\exists g \in G$  \; s.t. $\; g\cdot x=y$ \; (transitivity),

\medskip

4. $\int_M (\gamma(g) f) (x) d\mu(x)=\int_M f(g^{-1}\cdot x) d\mu(x)= \int_M  f(x) d\mu(x)\;\;$
$\forall g \in G\;$, $\forall f \in \bL^1$,

\smallskip

${}$ \; where $(\gamma(g) f )(x):= f(g^{-1}\cdot x)$.
Hence, $(\gamma(g))_{g\in G}$ is a group of isometry of $\bL^1$.

\begin{definition}
A continuous  real-valued kernel $K(x,y)$ on $M\times M$ is said to be $G$-invariant if
$$
K(g\cdot x, g\cdot y)= K(x,y), \quad \forall g\in G, \; \forall x,y \in M.
$$
If $K$ is the operator on $\bL^2$ with kernel $K(x, y)$, then $K$ is called $G-$invariant if
$\gamma(g) K = K \gamma(g)$, $\forall g \in G$, that is,
$$
\int_M K(g^{-1}\cdot x,y) f(y) d\mu(y) = \int_M K(x,y) f(g^{-1}\cdot y) d\mu(y),
\quad \forall  f \in \bL^2.
$$
\end{definition}

\begin{remark}\label{GM}
$(a)$ Assume that $K(x,y)$ is a continuous  $G$-invariant kernel, then

$(i)$ $ K(x,x)=K(g\cdot x, g\cdot x)$ and hence $K(x,x) \equiv   |M|^{-1}\Tr(K)$, and

\medskip

$(ii)$
$$
\int_M K(x,y) d\mu(y)= \int_M K(x, g\cdot y) d\mu(y)= \int_M K(g^{-1}\cdot x, y) d\mu(y),\quad \forall g\in G,
$$
and hence $\ONE:=\ONE_M$ is an eigenfunction of $K$, that is,
$\int_M K(x,y) \ONE(y) d\mu(y)= C \ONE(x)$.

\medskip

$(b)$ Suppose $K(x,y)$ is a continuous positive $G$-invariant kernel, then
$$
\psi_K(x,y) := K(x,x)+ K(y,y)-2K(x,y)= 2(C-K(x,y))
= 2(|M|^{-1}\Tr(K)-K(x,y)),
$$
clearly $\psi(x, y)$ is $G$-invariant and by (\ref{k1}), \;
$\tilde{K}(x,y) = K(x,y) + |M|^{-1}(\Tr(K)- 2C')$.

\medskip

$(c)$ Suppose $\psi(x,y)$ is a $G-$invariant N.D. kernel and
consider the associated  P.D. kernel $\tilde{K}$, defined as in \eref{def.tildeK}.
Then $\tilde{K}(x,y)$ is $G$-invariant, and
$$
x\mapsto  \frac 1{|M|}\int_M  \psi(x,u) d\mu(u) \equiv C_0  \;
\hbox{and }\;  \tilde{K}(x,y)= C_0- \frac 12\psi(x,y).
$$
Thus, in this framework there is one-to-one correspondence up to a constant between invariant P.D. and N.D. kernels.
\end{remark}

\subsection{Composition of operators}\label{comp}

Let $K(x,y)$ and $H(x,y)$ be two continuous kernels on $M\times M$ as above, and let $K$ and $H$ be the associate operators.
The operator $K\circ H$ is also a kernel operator with kernel $K\circ H(x,y)$:
$$
K\circ H(x,y)= \int_M K(x,u) H(u,y) d\mu(u).
$$
Observe that:
\begin{enumerate}
\item
If $K(x,y)= K(y,x) , H(x,y)= H(y,x)$ then
$$
K\circ H(x,y)= \int_M K(x,u) H(u,y) d\mu(u)=\int_M H(y,u) K(u,x)  d\mu(u)=
H\circ K(y,x).
$$
\item
If $K(x,y)$ and $H(x,y)$ are $G-$invariant, then so is $ K\circ H$. Indeed,
\begin{align*}
K\circ H(g\cdot x,g\cdot y)&= \int_M K(g\cdot x,u) H(u,g\cdot y) d\mu(u)
=\int_M K(g\cdot x,g\cdot u) H(g\cdot u, g\cdot y) d\mu(u)\\
&=\int_M K(x,u) H(u,y) d\mu(u)=K\circ H(x,y).
\end{align*}
\end{enumerate}

\subsection{Group action and metric}

Assume that we are in the setting of a Dirichlet space defined through
a  non-negative self-adjoint operator on $\bL^2(M,\mu)$ just as in \S\ref{subsec:background}.
Suppose  now that,
$$
\gamma(g) A= A\gamma(g), \quad \forall g \in G
$$
or equivalently
$$
\gamma(g) P_t= P_t \gamma(g), \quad \forall t>0, \;  \forall g \in G,
$$
i.e.  $\forall t>0, \; p_t(x,y)$ is $G-$invariant.
Clearly $\Gamma(f_1,f_2)$ is also $G$-invariant:
$\Gamma(f_1,f_2)=\Gamma(\gamma(g)f_1,\gamma(g)f_2)$
and the associate metric $\rho(x,y)$ is $G$-invariant:
$$
\rho(g\cdot x, g\cdot y)=\rho(x,y),\quad \forall g \in G.
$$

\begin{definition}
In the current framework, $(M, \mu, A, \rho, G)$ is said to be a two point homogeneous space if
$$
\forall x,y,x',y' \in M \;\; \hbox{s.t.}\;\;\rho(x,y)= \rho(x',y'), \; \exists g \in G\;\;
\hbox{s.t.} \;\; g\cdot x=x', \; g\cdot y=y'.
$$

In particular, $\forall (x,y) \in M \times M, \; \exists g \in G \; \hbox{s.t.} \; g\cdot x=y, \; g\cdot y=x$.
\end{definition}

\begin{theorem}\label{thm:two-point-hsp}

Let $(M, \mu, A, \rho, G)$ be a compact two point homogeneous space. Then we have:
\begin{enumerate}
\item
Any $G$-invariant continuous kernel $K(x,y)$ is symmetric.

\item
If $K(x,y)$ and $H(x,y)$ are two $G$-invariant continuous kernels, then $K\circ H= H\circ K$.

In particular, if $K(x,y)$ is a $G$-invariant continuous kernel, then $KA=AK$.

\item
Any  $G$-invariant real-valued continuous kernel $K(x,y)$ depends only on the distance $\rho(x,y)$, that is,
there exist a continuous function $k:\bR \mapsto \bR$,
such that
$$ K(x,y)= k(\rho(x,y)), \quad \forall x, y\in M.$$
\end{enumerate}
\end{theorem}

This theorem is a straightforward consequence of the observations from \S\ref{comp}
and the definition of two point homogeneous spaces.

\medskip

Let now $M$ be a compact Riemannian manifold and assume that $A:= -\Delta_M$ is the Laplacian on $M$,
$\rho$ is the Riemannian metric, and $\mu$ is the Riemannian measure.
Also, assume that there exists a compact Lie group $G$ of isometry on $M$ such that
$(M, \mu, -\Delta_M, \rho, G)$ is a compact two point homogeneous space.
 For the link with Gaussian processes see: \cite{AB}, \cite{G}.

Let $0\le \lambda_1 < \lambda_2 < \cdots$ be the spectrum of $-\Delta_M$.
Then the eigenspaces $\cH_{\lambda_k}:=\Ker (\Delta_M +\lambda_k \Id)$ are finite dimensional and
$$
\bL^2(M, \mu)= \bigoplus_{k\ge 1} \cH_{\lambda_k}.
$$
Let  $P_{\cH_{\lambda_k}}(x,y)$ be the kernel of the orthogonal projector onto $\cH_{\lambda_k}$.
Then if $K(x,y)$ is a $G$-invariant positive definite kernel
we have the following decomposition of $K(x,y)$, which follows from Bochner-Godement theorem (\cite{Far-Har}, \cite{HEL}):
$$
K(x,y)= \sum_{k\geq 0} \nu_k P_{\cH_{\lambda_k}}(x,y), \quad \nu_k \geq 0.
$$

\section{Brownian motion on the sphere}\label{sec:sphere}

In this section we apply our main result (Theorem~\ref{MAIN}) to a Gaussian process parametrized by
the unit sphere $\bS^d$ in $\bR^{d+1}$.
This is a Riemannian manifold and a compact two point homogeneous space.
More explicitly,
$$
G=SO(d+1), \; H= SO(d), \; G/H = \bS^d.
$$
The geodesic distance $\rho$ on $\bS^d$ is given by
$$
\rho(\xi, \eta) = \arccos \langle \xi, \eta \rangle,
$$
where $\langle \xi, \eta \rangle$ is the inner product of $\xi, \eta\in \bR^{d+1}$.
Clearly,
$$
\forall \xi, \eta \in \bS^d, \; \forall g \in G,\;\;
\rho(g\cdot \xi, g\cdot \eta)= \rho(\xi, \eta),
\quad\hbox{and}\quad
\forall \xi, \eta \in \bS^d,\; \exists g\in G \;\;\hbox{s.t.}\;\; g\cdot \xi=\eta.
$$
Thus $G$ acts isometrically and transitively on $\bS^d$.
Furthermore,
$$\hbox{
$\forall \xi, \eta, \xi', \eta' \in \bS^d$ \;\;s.t. \;\;$\rho(\xi, \eta)=\rho(\xi', \eta')$,
there exists $g \in G$ \;\;s.t.\;\; $g\cdot \xi= \xi'$ and $g\cdot \eta = \eta'$.
}
$$
Therefore, $\bS^d$ is a compact two point homogeneous space.

Let $-\Delta_{\bS^d}$ be the (positive) Laplace-Beltrami operator on $\bS^d$.
As is well known the eigenspaces of $-\Delta_{\bS^d}$ are the spaces of spherical harmonics, defined by
$$
\cH_{\lambda_k} := \Ker (\Delta_{\bS^d} +\lambda_k I_d),
\quad \lambda_k := k(k+d-1) = k(k+2\nu), \;\; k\ge 0 \quad \nu:= \frac{d-1}2.
$$
One has $\bL^2(\bS^d)= \bigoplus_{k\ge 0} \cH_{\lambda_k}$
and the kernel of the orthogonal projector $P_{\cH_{\lambda_k}}$ onto $\cH_{\lambda_k}$ is given by
$$
P_{\cH_{\lambda_k}}(\xi,\eta)
= L_k^d( \langle \xi, \eta \rangle), \quad L_k^d(x):=|\bS^d|^{-1}\Big(1+ \frac k\nu \Big)C_k^\nu(x).
$$
Here $C^\nu_k(x)$, $k\ge 0$, are the Gegenbauer polynomials defined on $[-1, 1]$ by the generating function
$$
\frac 1{ (1-2xr+r^2)^\nu} = \sum_{k\ge 0} r^k C^\nu_k(x).
$$
Therefore,
$$
-\Delta_{\bS^d}f = \sum_{k\ge 0} k(k+2\nu) P_{\cH_{\lambda_k}} f
$$
and the invariant continuous positive definite functions on $\bS^d$ are of the form
$$
K(\xi,\eta)= \sum_k \nu_k  L_k^d( \langle \xi, \eta \rangle)
= \sum_k \nu_k  L_k^d( \cos \rho(\xi, \eta)),
$$
where
$$
\sum_k \nu_k  L_k^d(1) =  \sum_k \nu_k  L_k^d( \langle \xi, \xi \rangle)  <\infty.
$$
Note that
$$
L_k^\nu(1) |\bS^d| = \int_{\bS^d}  L_k^\nu( \langle \xi, \xi \rangle)d\mu(\xi)
= \dim(\cH_{\lambda_k}(\bS^d))= \binom{k+d}{d}- \binom{k-2+d}{d}\sim k^{d-1}.
$$
Let
$$
W^\nu_k(x):= \frac{L_k^\nu(x)}{L_k^\nu(1)}
= \frac{C_k^\nu(x)}{C_k^\nu(1)}.
\quad \hbox{Clearly,} \;\;
W^\nu_k(1)=\sup_{x\in [-1,1]} |W^\nu_k(x)|=1.
$$
Then (see \cite{B})
$$
\lim_{\nu \mapsto 0} \frac{C_k^\nu (x)}{C_k^\nu(1)} = T_k(x) \;\; (= W^{0}_k(x) \; \hbox{by convention} ),
$$
$$
\lim_{\nu \mapsto \infty} \frac{C_k^\nu(x)}{C_k^\nu(1)}  = x^k \;\; (= W^\infty_k(x) \; \hbox{by convention} ).
$$
Here $T_k $ is the Chebyshev polynomial of first kind ($T_k(\cos \theta)= \cos k\theta$).
The invariant continuous positive definite functions on $\bS^d$ are of the form
$$
K^\nu(\xi,\eta)= \sum_{k\ge 0} a^\nu_k  W_k^\nu( \langle \xi, \eta \rangle)
= \sum_{k\ge 0} a_k^\nu  W_k^\nu(\cos \rho(\xi, \eta)), \quad a^\nu_k\geq 0, \quad  \sum_k a_k^\nu   <\infty.
$$
Clearly,
\begin{equation}\label{W}
 \sum_k a_k^\nu  W_k^\nu( \cos \rho(\xi, \eta))
 =  \sum_k  \frac{a_k^\nu}{ L_k^\nu(1)}L_k^\nu( \cos \rho(\xi, \eta)),
\quad  L_k^\nu(1) \sim k^{d-1}.
 \end{equation}
Therefore,
$$
\nu_k=|\bS^d|\frac{a_k^\nu}{\dim( \cH_{\lambda_k})} =O\Big(\frac{a_k^\nu}{ k^{d-1}}\Big).
$$

The following {\bf Schoenberg-Bingham result} (see e.g. \cite{B}) plays a key role here:
{\em
If $f $ is a continuous function defined on $[-1,1]$,
then $f(\langle \xi, \eta \rangle)$ is a positive definite function on $\bS^d$
and invariant with respect to $SO(d+1)$ for all $d \in \bN$
if and only if
$$
f(x) = \sum_{n \geq 0} a_n x^n, \quad \hbox{where}\quad a_n\ge 0
\quad\hbox{and}\quad \sum_{n\ge 0} a_n = f(1) <\infty.
$$
}

Therefore, for such a function $f$
$$
f(x)= \sum_{k \geq 0} a_k^\nu W^\nu_k(x), \quad a_k^\nu\ge 0,
\quad\hbox{and} \quad  \sum_{k \geq 0}  a_k^\nu =\sum_{k\ge 0} a_k = f(1),
$$
and hence
$$
f(\langle  \xi, \eta  \rangle)=  \sum_{k \geq 0} a_k^\nu  W_k^\nu( \langle \xi, \eta  \rangle)
= \sum_{k \ge 0} \frac{ a_k^\nu}{L_k^\nu(1)}  L_k^\nu( \langle \xi, \eta  \rangle)
= f(\cos \rho(\xi, \eta)).
$$

\subsection{Fractional Brownian process on the sphere}

\begin{theorem}\label{thm:sphere-1}
For any $0<\alpha\le 1$ the function
$$
\psi(\xi, \eta)= \rho(\xi, \eta)^\alpha, \quad \xi, \eta\in \bS^d,
$$
is negative definite,
and the associated Gaussian process has almost everywhere regularity
$B^\gamma_{\infty,1}$, $\gamma <\frac \alpha 2$.
\end{theorem}

\noindent
{\bf Proof.}
Consider first the case when $\alpha=1$ ({\em Brownian process}).
We will show that for some constant $C>0$ the function $C- \rho(\xi, \eta)$
is an invariant positive definite function.
To this end, by Schoenberg-Bingham result
we have to prove that there exists a function
$$
f(x)= \sum a_n x^n, \quad\hbox{with}\;\; a_n\ge 0, \;\; \sum_{n\ge 0} a_n <\infty,
$$
such that
$
f(\langle  \xi, \eta  \rangle)= f(\cos \rho(\xi, \eta))= C-  \rho(\xi, \eta).
$
Luckily the function $\frac \pi2- \arccos x$ does the job.
Indeed, it is easy to see that
$$
f(x):= \frac \pi2- \arccos x =\arcsin x
= \sum_{j\ge 0} \frac{ (\frac 12)_j(\frac 12)_j}{j!(\frac 32)_j}x^{2j+1}
\quad \hbox{and}\quad
\sum_{j\geq 0} \frac{ (\frac 12)_j(\frac 12)_j }{j! (\frac 32)_j } = \frac \pi2\;  \; \hbox{(Gauss)}.
$$
Here we use the standard notation
$(a)_j := a(a+1)\cdots(a+j-1)=\Gamma(a+j)/\Gamma(a)$.
Therefore,
$$
f(\langle  \xi, \eta  \rangle)
= \frac \pi2- \arccos \langle  \xi, \eta  \rangle
= \frac{\pi}2- \rho(\xi, \eta).
$$
Clearly,
$
|f(\langle  \xi, \eta  \rangle)-f(\langle  \xi, \eta'  \rangle)| \le \rho(\eta, \eta')
$
and by Theorem~\ref{MAIN} the associated Gaussian process $(Z^d_\xi(\omega))_{\xi \in \bS^d}$
is almost surely in $B^{s}_{\infty, 1}(\bS^d)$ (hence in $\Lip s$)
for $0< s<\frac 12$.
Furthermore,
$$
\bE(Z^d_\xi-Z^d_\eta)^2 = 2f(1)- 2f(\langle  \xi, \eta  \rangle)= 2 \rho(\xi, \eta).
$$

Consider now the general case: $0<\alpha\le 1$ ({\em Fractional Brownian process}).
From above it follows that
$\psi(\xi, \eta):= \rho(\xi, \eta)$ is an invariant negative definite kernel.
Then the general theory of negative definite kernels yields that for any $0<\alpha\le 1$ the kernel
$\psi_\alpha(\xi, \eta)= \rho(\xi, \eta)^\alpha$
is invariant and negative definite.
Therefore, for a sufficiently large constant $C>0$,
$$
K(\xi, \eta) = C- \frac 12 \rho(\xi, \eta)^\alpha
$$
is an invariant positive definite kernel.
On the other hand,
$$
|K(\xi, \eta) -K(\xi, \eta') |
=  \frac 12 |\rho(\xi, \eta)^\alpha - \rho(\xi, \eta')^\alpha|
\le \frac 12 \rho(\eta', \eta))^\alpha.
$$
By Theorem~\ref{MAIN} it follows that the associated Gaussian process $(Z^d_\xi(\omega))_{\xi \in \bS^d}$
is almost surely in  $B^\gamma_{\infty,1}$, $\gamma <\frac \alpha 2$,
and hence in $\Lip s$, $s<\frac \alpha 2$, and the proof is complete.
\qed

\begin{remark}\label{rem:csensov}
From the definition of the process, we have
$$
\bE(Z_\xi^\alpha -Z_\eta^\alpha)^2 = \rho(\xi, \eta)^\alpha.
$$
This directly connects to the regularity proofs of such a process using generalization of Kolmogorov-Csensov inequalities.
 See for instance \cite{AL} and \cite{LSch}.
\end{remark}

\subsection{Regularity of Gaussian processes on the sphere: General result}

\begin{theorem}
Let
$$
f(x) = \sum_{n\ge 0} \frac{A_n}{n!} x^n,
\quad\hbox{where}\quad  A_n \geq 0, \;\;\hbox{and}\;\;
\frac{A_n}{n!} =O\big(\frac 1{n^{1+\alpha}}\big), \quad \alpha >0.
$$
Then
$$
K(\xi, \eta):= f(\cos \langle \xi, \eta \rangle), \;\; \xi, \eta \in \bS^d, \;\; d\ge 1,
$$
is an invariant positive definite function,
and the associated Gaussian process $(Z^d_\xi(\omega))_{\xi \in \bS^d}$ is almost surely in
$B^\gamma_{\infty, 1}$ for $\gamma < \alpha$.
\end{theorem}

\noindent
{\bf Proof.}
By Corollary~\ref{dec}, it suffices to show that $f(x)$ can be represented in the following form (see (\ref{W})):
$$
f(x)= \sum_j B_j W^\nu_j(x), \quad   0 \le B_j = O\Big(\frac 1{j^{1+2\alpha}}\Big),
\quad\hbox{implying}\;\;
\nu_j =O \Big(\frac 1{j^{d+2\alpha}}\Big)= O( \sqrt{\lambda_j} )^{2\alpha +d}.
$$
By lemma~1 in \cite{B} and the obvious identity $\Gamma(x+n)=(x)_n\Gamma(x)$
we obtain the representation
$$
x^n =  \frac{n!}{2^n}  \sum_{0\le 2k\le n}
 \frac {n-2k+\nu}{k!(\nu)_{n-k +1}}\frac{(2\nu)_{n-2k}}{(n- 2k)!}W^\nu_{n-2k}(x).
$$
Substituting this in the definition of $f(x)$ we obtain
\begin{align*}
f(x)= \sum_{n\ge 0}\frac{A_n}{n!}x^n
&= \sum_{n\ge 0} \frac{A_n}{2^n}\sum_{0\le 2k\le n}
\frac{n-2k+\nu}{k!(\nu)_{n-k +1}}\frac{(2\nu)_{n-2k}}{(n- 2k)!}W^\nu_{n-2k}(x) \quad\quad(j=n-2k)
\\
&= \sum_{j\geq 0} \frac{(j+\nu)(2\nu)_{j}}{j!}  W^\nu_{j}(x)
\sum_{n-2k =j} \frac{A_n}{2^nk!(\nu)_{n-k +1}}
\\
&= \sum_{j\geq 0} \frac{(j+\nu)(2\nu)_{j}}{j!} W^\nu_{j}(x) \frac 1{2^j}
\sum_{k\ge 0} \frac{A_{j+2k}}{2^{2k}k!(\nu)_{j+k +1}}
 =: \sum_{j\ge 0} B_j W^\nu_{j}(x),
\end{align*}
where
\begin{align*}
B_j &:= \frac{(j+\nu) (2\nu)_{j}}{j!2^j} \sum_{k\ge 0} \frac{A_{j+2k}}{2^{2k}k!(\nu)_{j+k +1}}
\\
&= \frac{(j+\nu)(2\nu)_{j}}{j!2^j(\nu)_{j+1}}
\sum_{k\ge 0} \frac{A_{j+2k}}{2^{2k}k!(\nu+j+1)_{k}}
\\
&= \frac{(2\nu)_j}{2^j j!(\nu)_j}
\sum_{k\geq 0} \frac{A_{j+2k}}{2^{2k}k!(\nu+j+1)_{k}}.
\end{align*}
However, for $n>\alpha$ we have
\begin{align*}
\frac{c_1(\alpha)}{n^{1+\alpha}} \le \frac{\Gamma(n-\alpha)}{n!} \le \frac{c_2(\alpha)}{n^{1+\alpha}}
\quad\hbox{and hence}\quad
\frac{A_n}{n!} = O\Big( \frac 1{n^{1+\alpha}} \Big) \; \Longleftrightarrow \;
A_n = O(\Gamma(n-\alpha)).
\end{align*}
We use this to obtain for $j>\alpha$ (with $c=c(\alpha)$)
\begin{align*}
\sum_{k\geq 0} \frac{A_{j+2k}}{2^{2k}k!(\nu+j+1)_{k}}
&\le c\sum_{k\geq 0} \frac{\Gamma(j+2k-\alpha)}{2^{2k}k!(\nu+j+1)_{k}}
\\
&= c\Gamma (j-\alpha)
 \sum_{k\geq 0} \frac{\Gamma( j+2k- \alpha)}{ \Gamma (j-\alpha)}
  \frac{1}{2^{2k}k! (\nu+j+1)_{k}}
\\
&= c\Gamma (j-\alpha)
 \sum_{k\geq 0}
  \frac{(j- \alpha)_{2k}}{2^{2k}}  \frac 1{ k! (\nu+j+1)_{k}}
\\
&= c\Gamma (j-\alpha)
 \sum_{k\geq 0}  \Big(\frac{j- \alpha}2\Big)_{k} \Big(\frac{j- \alpha +1}2\Big)_{k}
   \frac 1{ k! (\nu +j+1)_{k }},
\end{align*}
where we used the Legendre duplication formula (see e.g. \cite{AAR}):
$$
\frac{(b)_{2k}}{2^{2k}} = \frac{\Gamma(b+2k)}{ 2^{2k}\Gamma (b)}
=\Big(\frac b2\Big)_k \Big(\frac{b+1}2\Big)_k.
$$
By the Gaussian identity (see e.g. \cite[Theorem 2.2.2]{AAR})
\begin{align*}
\sum_{k\geq 0}\Big(\frac{j- \alpha}2\Big)_{k} \Big(\frac{j- \alpha +1}2\Big)_{k} \frac 1{ k! (\nu +j+1)_{k}}
&= \frac{ \Gamma(\nu+j+1)  \Gamma(\nu+j+1-\frac{j- \alpha}2 -\frac{j- \alpha +1}2  )  }{ \Gamma(\nu+j+1-\frac{j- \alpha}2)  \Gamma(\nu+j+1-\frac{j- \alpha +1}2)}
\\
&= \frac{ \Gamma(\nu+j+1)  \Gamma(\nu+\frac 12 + \alpha   )  }
    { \Gamma(\nu+  \frac j2+1 +  \frac \alpha2)  \Gamma(\nu+ \frac j2+ \frac 12 +  \frac \alpha2)}
\end{align*}
and hence
$$
B_j \le c\frac{(2\nu)_j}{j!2^j(\nu)_j}
\frac{\Gamma (j-\alpha)\Gamma(\nu+j+1)\Gamma(\nu+\frac 12 + \alpha)}
{\Gamma(\nu+ \frac j2+1 + \frac \alpha2)\Gamma(\nu+ \frac j2+ \frac 12 +  \frac \alpha2)}.
$$
Applying again the Legendre duplication formula, we get
$$
\Gamma\Big(\frac 12\Big) \Gamma (2\nu+j+1+\alpha)
= \Gamma\Big(\nu+  \frac j2+1 +  \frac \alpha2\Big)
\Gamma\Big(\nu+ \frac j2+ \frac 12 +  \frac \alpha2\Big) 2^{2\nu+j+\alpha}.
$$
We use this above to obtain for $j\ge 2\alpha$
\begin{align*}
B_j &\le c\frac{(2\nu)_j}{j!(\nu)_j}
\frac{\Gamma (j-\alpha)\Gamma(\nu+j+1)\Gamma(\nu+\frac 12 + \alpha)}
{\Gamma(\frac 12)\Gamma(2\nu+j+1+\alpha)2^{-2\nu-\alpha}}
\\
& = c\frac{\Gamma(2\nu+j)\Gamma(\nu)}{\Gamma(j+1)\Gamma(2\nu)\Gamma(\nu+j)}
\frac{\Gamma (j-\alpha)\Gamma(\nu+j+1)\Gamma(\nu+\frac 12 + \alpha)}
{\Gamma(\frac 12)\Gamma(2\nu+j+1+\alpha)2^{-2\nu-\alpha}}
\\
&= c2^{\alpha+1}(j+\nu)
\frac{\Gamma(\nu+\frac 12 + \alpha)}{\Gamma(\nu +\frac 12)}
\frac{\Gamma (j-\alpha) }{\Gamma(j-\alpha +1+ \alpha)}
\frac{\Gamma(2\nu+j)}{\Gamma(2\nu+j+1+\alpha)}
\\
&\le c(j+\nu) \frac 1{(j-\alpha)^{1+\alpha}}  \frac 1{(2\nu +j)^{1+\alpha}}
\le \frac{c}{j^{1+2\alpha }}.
\end{align*}
Here we used once again the the Legendre duplication formula.
It is easy to show that $B_j \le c(\alpha)$, if $j<2\alpha$.
Therefore,
$B_j = O\big(\frac 1{j^{1+2\alpha}}\big)$
and this completes the proof.\qed

\begin{corollary}
Let $a> 0$, $b>0$, $c> a+b$, $\alpha= c-a-b$, and let
$$
F_{a,b;c}(x):= \sum_n \frac{(a)_n (b)_n}{({c})_n} \frac{x^n}{n!}.
$$
Then
$F_{a,b;c}(\langle \xi, \eta \rangle)$
is an invariant positive definite function on the sphere $\bS^d$
and the associated Gaussian process has regularity
$B^\gamma_{\infty,1}, \; \gamma < \alpha ,$ almost everywhere.
\end{corollary}


\begin{thebibliography}{99}

\bibitem{Adler1990}
Adler, R.J.: An introduction to continuity, extrema, and related topics for general Gaussian processes.
Institute of Mathematical Statistics Lecture Notes—Monograph Series, 12.
Institute of Mathematical Statistics, Hayward, CA (1990)

\bibitem{Adler}
Adler, R.J., Taylor J.E.: Random fields and geometry.
Springer Monographs in Mathematics. Springer, New York (2007)

\bibitem{AL}
Andreev, R., Lang, A.:
Kolmogorov-Chentsov theorem and differentiability of random field on manifolds.
Potential Anal. {\bf 41}(3), 761–-769 (2014)

\bibitem{AAR}
Andrew, G., Askey, R., Roy, R.:
Special functions.
Encyclopedia of Mathematics and its Applications, 71. Cambridge University Press, Cambridge (1999)

\bibitem{ARZ}
Aronszajn, N.:
Theory of reproducing kernels. Trans. Amer. Math. Soc. {\bf 68}(3), 337–-404 (1950)

\bibitem{AB}
Askey, R., Bingham, N. H.: Gaussian processes on compact symmetric spaces.
Z. Wahrscheinlichkeitstheorie und Verw. Gebiete {\bf 37}(2), 127-–143 (1976/77)

%
\bibitem{BCR}
Berg, C., Christensen, J.P.R., Ressel, P.:
Harmonic analysis on semigroups.
Theory of positive definite and related functions. Graduate Texts in Mathematics, 100.
Springer-Verlag, New York (1984)

\bibitem{B}
Bingham, N.H.:
Positive definite functions on spheres. Proc. Cambridge Philos. Soc. {\bf 73}, 145–-156 (1973)

\bibitem{Bochner}
Bochner, S.:
Harmonic analysis and the theory of probability.
University of California Press, Berkeley and Los Angeles (1955)

\bibitem{Bogachev}
Bogachev, V.I.: Gaussian Measures.
Mathematical Surveys and Monographs, 62. American Mathematical Society, Providence, RI (1998)

\bibitem{BH}
Bouleau, N., Hirsch, F.:
Dirichlet forms and analysis on Wiener space.
de Gruyter Studies in Mathematics, 14. Walter de Gruyter \& Co., Berlin (1991)

\bibitem{CKR}
Ciesielski, Z. ,Kerkyacharian, G., Roynette, B.:
Quelques espaces fonctionnels associes a des processus Gaussiens.
Studia Math. {\bf{107}}(2), 171--204 (1993)

\bibitem{CW}
Coifman, R., Weiss, G.:
Analyse harmonique non-commutative sur certains espaces homog\`{e}nes.
Lecture Notes in Mathematics, Vol. 242. Springer-Verlag, Berlin-New York (1971)

\bibitem{CKP}
Coulhon, T., Kerkyacharian, G., Petrushev, P.:
Heat kernel generated frames in the setting of Dirichlet spaces,
J. Fourier Anal. Appl. {\bf{18}}(5), 995--1066 (2012)

%
\bibitem{CUS}
Cucker, F., Smale, S.:
On the mathematical foundations of learning.
Bull. Amer. Math. Soc. (N.S.) {\bf 39}(1), 1-–49 (2002)

%
\bibitem{Davies1}
Davies, E. B.:
Linear operators and their spectra.
Cambridge Studies in Advanced Mathematics, 106.
Cambridge University Press, Cambridge (2007)

\bibitem{Far-Har}
Faraut, J., Harzallah, K.:
Distances hilbertiennes invariantes sur un espace homog\`{e}ne.
Ann. Inst. Fourier (Grenoble) {\bf 24}(3), 171–-217 (1974)

\bibitem{Fernique}
Fernique, X.:
Regularit\'{e} des trajectoires des fonctions al\'{e}atoires gaussiennes.
\'{E}cole d'\'{E}t\'{e} de Probabilit\'{e}s de Saint-Flour, IV-1974, pp. 1–-96.
Lecture Notes in Math., Vol. 480, Springer, Berlin (1975)

%
%
%
\bibitem{FUKU}
Fukushima, M., Oshima, Y., Takeda, M.:
Dirichlet forms and symmetric Markov processes.
Second revised and extended edition.
de Gruyter Studies in Mathematics, 19. Walter de Gruyter \& Co., Berlin (2011)

\bibitem{G}
Gangolli, R.:
Positive definite kernels on homogeneous spaces and certain stochastic processes
related to L\'{e}vy's Brownian motion of several parameters.
Ann. Inst. H. Poincar\'{e} Sect. B (N.S.) {\bf 3}(2), 121-–226 (1967)

\bibitem{Grigorian}
Grigor'yan, A.:
Heat kernel and analysis on manifolds.
AMS/IP Studies in Advanced Mathematics, 47.
American Mathematical Society, Providence, RI; International Press, Boston, MA (2009)


\bibitem{HEL}
Helgason, S.:
Differential geometry and symmetric spaces.
Pure and Applied Mathematics, Vol. XII. Academic Press, New York-London (1962)

\bibitem{GKPP}
Kerkyacharian, G., Petrushev, P.:
Heat kernel based decomposition of spaces of distributions in the framework of Dirichlet spaces.
Trans. Amer. Math. Soc. {\bf 367}(1), 121-–189 (2015)

\bibitem{VV2}
Knapik, B. T., van der Vaart, A. W.,  van Zanten, J. H.:
Bayesian inverse problems with Gaussian priors.
Ann. Statist. {\bf 39}(5), 2626–-2657 (2011)

\bibitem{LSch}
Lang, A., Schwab, C.:
Isotropic Gaussian random field on the sphere:
regularity, fast simulation, and stochastic partial
differential equation.
arxiv:1305.1170v2. 16 May 2014.

\bibitem{Ledoux}
Ledoux, M.:
Isoperimetry and Gaussian analysis.
Lectures on probability theory and statistics (Saint-Flour, 1994), 165–-294,
Lecture Notes in Math., 1648, Springer, Berlin (1996)

\bibitem{Led}
Ledoux, M.:
The concentration of measure phenomenon.
Mathematical Surveys and Monographs, 89.
American Mathematical Society, Providence, RI (2001)

\bibitem{LedTala}
Ledoux, M., Talagrand, M.:
Probability in Banach spaces. Isoperimetry and processes.
Ergebnisse der Mathematik und ihrer Grenzgebiete (3), 23.
Springer-Verlag, Berlin (1991)

\bibitem{Li}
Li, W.V., Shao, Q.M.:
Gaussian processes: inequalities, small ball probabilities and applications.
Stochastic processes: theory and methods, 533--597, Handbook of Statistcs 19,
North-Holland, Amsterdam (2001)

\bibitem{Lifshits}
Lifshits, M.:
Lectures on Gaussian processes. Springer Briefs in Mathematics. Springer, Heidelberg (2012)

\bibitem{Lus-Pages}
Luschgy, H., Pages, G.:
Expansions for Gaussian processes and Parseval frames.
Electron. J. Probab. {\bf 14}(42), 1198–-1221 (2009)

\bibitem{Rosen}
Marcus, M.B., Rosen, J.:
Markov processes, Gaussian processes, and local times.
Cambridge Studies in Advanced Mathematics, 100. Cambridge University Press, Cambridge (2006)

%
%
%
\bibitem{Peetre}
Peetre, J.:
New thoughs on Besov spaces,
Duke University (1976)

%


\bibitem{Rasmussen}
Rasmussen, C. E., Williams, Christopher K. I.:
Gaussian processes for machine learning.
Adaptive Computation and Machine Learning. MIT Press, Cambridge, MA (2006)


%
%

\bibitem{Schoenberg1}
Schoenberg, I.J.:
Metric spaces and positive definite functions.
Trans. Amer. Math. Soc. {\bf{44}}(3), 522--536 (1938)

\bibitem{Seeger}
Seeger, M.:
Gaussian processes for machine learning.
International Journal of Neural Systems, {\bf 14}(02), 69--106 (2004)

\bibitem{SteinChrist}
Steinwart, I., Christmann, A.:
Support vector machines.
Information Science and Statistics. Springer, New York (2008)

\bibitem{Sturm0}
Sturm, K.T.:
Analysis on local Dirichlet spaces I. Recurrence, conservatisness and $L_p$-Liouville properties,
J. Reine Angew. Math. {\bf{456}}, 173--196 (1994)

\bibitem{Sturm1}
Sturm, K.T.:
Analysis on local Dirichlet spaces II.
Upper Gaussian estimates for the fundamental solutions of parabolic equations,
Osaka J. Math. {\bf{32}}(2), 275--312 (1995)

\bibitem{Sturm}
Sturm, K.T.:
Analysis on local Dirichlet spaces III. The parabolic Harnack inequality,
J. Math. Pures Appl. {\bf{75}}(3), 273--297 (1998)

\bibitem{Sturm3}
Sturm, K.T.:
The geometric aspect of Dirichlet forms.
New directions in Dirichlet forms, 233--277,
AMS/IP Stud. Adv. Math., {\bf{8}}, Amer. Math. Soc., Providence, RI(1998)

\bibitem{Talagrand}
Talagrand, M.:
Regularity of Gaussian processes.
Acta Math. {\bf{159}}(1-2), 99--149 (1987)

\bibitem{Ta}
Talagrand, M.:
Mean field models for spin glasses,
Springer-Verlag, Berlin (2011)

\bibitem{Talagrand3}
Talagrand, M.:
Upper and lower bounds for stochastic processes.
Springer, Heidelberg (2014)

%
\bibitem{VV1}
van der Vaart, A., van Zanten, H.:
Bayesian inference with rescaled Gaussian process priors.
Electron. J. Stat. {\bf 1}, 433--448 (2007)


\end{thebibliography}
\end{document}